\documentclass[a4paper,11pt]{amsart}
\usepackage{amssymb,amscd}
\theoremstyle{plain}

\newtheorem{theo}{Theorem}[section]
\newtheorem{prop}[theo]{Proposition}
\newtheorem{lemm}[theo]{Lemma}
\newtheorem{coro}[theo]{Corollary}

\theoremstyle{definition}
\newtheorem*{defi}{Definition}
\newtheorem{exam}[theo]{Example}

\theoremstyle{remark}
\newtheorem*{rema}{Remark}

\numberwithin{equation}{section}

\newcommand{\mb}[1]{{\textbf {\textit#1}}}

\newcommand{\field}[1]{\mathbb{#1}}

\newcommand{\R}{\field{R}}

\newcommand{\Z}{\field{Z}}
\renewcommand{\k}{\mathbf k}

\DeclareMathOperator{\lk}{lk}
\DeclareMathOperator{\rank}{rk}

\DeclareMathOperator{\Hom}{Hom}

\DeclareMathOperator{\Map}{Map}
\DeclareMathOperator{\st}{st}
\DeclareMathOperator{\Tor}{Tor}

\def\G{\Gamma}
\def\mG{\Gamma}

\def\mF{F}

\def\o{o}
\def\Htwo{H^2(BT)}
\def\HT{H^*_T}
\def\SV{H^*(BT)}

\newcommand{\ta}[1]{\tau_{{}_{\!#1}} }
\renewcommand{\v}[1]{ v_{{}_{\!#1}} }

\def\hatzero{\hat 0}
\def\P{\mathcal P}

\def\le{\leqslant}
\def\ge{\geqslant}

\newcommand{\n}[1]{\nu_{{}_{\!#1}} }
\def\til{\widetilde}
\renewcommand{\u}[1]{ u_{{}_{\!#1}} }

\begin{document}

\title{Torus graphs and simplicial posets}

\author{Hiroshi Maeda}
\author{Mikiya Masuda}
\author{Taras Panov}\thanks{The third author was supported by the
21 COE Programme at Osaka City University,
the Japanese Society for the Promotion of Science (grant
no.~P05296), and the Russian Foundation
for Basic Research (grant no.~04-01-00702).}
\address{Hakuryo Co., Ltd}
\address{Department of Mathematics, Osaka City University,
Sumiyoshi-ku, Osaka 558-8585, Japan}
\address{Department of Geometry and Topology, Faculty of Mathematics and Mechanics, Moscow State University,
Leninskiye Gory, Moscow 119992, Russia\newline \emph{and}
\newline Institute for Theoretical and Experimental Physics, Moscow
117259, Russia}
\email{maeda.hiroshi@hakuryo.com}
\email{masuda@sci.osaka-cu.ac.jp}
\email{tpanov@mech.math.msu.su}


\begin{abstract}
For several important classes of manifolds acted on by the torus,
the information about the action can be encoded combinatorially by
a regular $n$-valent graph with vector labels on its edges, which
we refer to as the \emph{torus graph}. By analogy with the
\emph{GKM-graphs}, we introduce the notion of \emph{equivariant
cohomology} of a torus graph, and show that it is isomorphic to
the \emph{face ring} of the associated \emph{simplicial poset}.
This extends a series of previous results on the equivariant
cohomology of torus manifolds. As a primary combinatorial
application, we show that a simplicial poset is
\emph{Cohen--Macaulay} if its face ring is Cohen--Macaulay. This
completes the algebraic characterisation of Cohen--Macaulay posets
initiated by Stanley. We also study \emph{blow-ups} of torus
graphs and manifolds from both the algebraic and the topological
points of view.
\end{abstract}

\maketitle

\section{Introduction}
The study of torus actions on manifolds is renowned for its close
connections with combinatorics and combinatorial geometry. Two
classes of actions are typical here; namely, the (smooth, compact)
\emph{algebraic toric varieties} and \emph{Hamiltonian} torus
actions on \emph{symplectic manifolds}. Both are very special
cases of a torus action $T^k\times M^{2n}\to M^{2n}$ on an
even-dimensional manifold; and the relation to combinatorics comes
from the study of the orbit poset and the orbit quotient. In the
former case the notion of \emph{fan}, which encodes both
combinatorial and geometric data, allows one to completely
translate algebraic geometry into combinatorics; in the latter
case important information about the Hamiltonian action is
contained in the \emph{moment polytope}.

During the last two decades both examples have developed into
several other classes of manifolds with torus action, mostly of
purely topological nature. These manifolds are neither algebraic
varieties nor symplectic manifolds in general, thereby enjoying a
larger flexibility for topological or combinatorial applications,
but still possess most of the important topological properties of
their algebraic or symplectic predecessors. The study of toric
varieties from a topological viewpoint led to the appearance of
\emph{(quasi)toric manifolds}~\cite{da-ja91},
\emph{multifans}~\cite{masu99} and \emph{torus
manifolds}~\cite{ha-ma03}. The latter carry an effective
half-dimensional torus action $T^n\times M^{2n}\to M^{2n}$ whose
fixed point set is non-empty.

The concept of a \emph{GKM-manifold} is closely related to
Hamiltonian torus actions. According to~\cite{gu-za99}, a compact
$2n$-dimensional manifold $M$ with an effective torus action
$T^k\times M\to M$ ($k\le n$) is called a GKM-manifold if the
fixed point set is finite, $M$ possesses an invariant almost
complex structure, and the weights of the tangential
$T^k$-representations at the fixed points are pairwise linearly
independent. These manifolds are named after Goresky, Kottwitz and
MacPherson, who studied them in~\cite{g-k-m98}. They showed that
the ``one-skeleton'' of such a manifold $M$, that is, the set of
points fixed by at least a codimension-one subgroup of $T^k$, has
the structure of a ``labelled'' graph $(\Gamma,\alpha)$, and that
the most important topological information about $M$, such as its
Betti numbers or equivariant cohomology ring, can be read directly
from this graph. These graphs have since become known as
\emph{GKM-graphs} (or \emph{moment graphs}); and their study has
been of independent combinatorial interest since the appearance of
Guillemin and Zara's paper~\cite{gu-za99}. The idea of associating
a labelled graph to a manifold with a circle action also featured
in Musin's work~\cite{musi80}.

Both GKM- and torus manifolds have become objects of study in the
emerging field of \emph{toric topology}, and linking these
important classes of torus actions together has been one of our
aims here. Our concept of a \emph{torus graph}, motivated by that
of a GKM-graph, allows us to translate the important topological
properties of torus manifolds into the language of combinatorics,
like in the case of GKM-manifolds. Therefore, the study of torus
graphs becomes our primary objective. A torus graph is a finite
$n$-valent graph $\Gamma$ (without loops, but with multiple edges
allowed) with an \emph{axial function} on the set $E(\Gamma)$ of
oriented edges taking values in $\Hom (T^n,S^1)=H^2(BT^n)$ and
satisfying certain compatibility conditions. These conditions
(described in Section~\ref{stogr}) are similar to those for
GKM-graphs, but not exactly the same. For the graphs coming from
torus manifolds the values of the axial function coincide with the
weights of the tangential $T^n$-representations at the fixed
points.

The notion of equivariant cohomology of a torus graph introduced
in Section~\ref{stogr} is same as that of a GKM-graph given
in~\cite{gu-za99} and~\cite{gu-za01}. However, unlike the case of
GKM-graphs, we have been able to completely describe the
equivariant cohomology ring of a torus graph in terms of
generators and relations, by applying the methods of our previous
work~\cite{ma-pa06} to the associated \emph{simplicial poset}.
Simplicial posets have already shown their importance in the
topological study of torus actions (see e.g.~\cite{ma-pa06}); they
also feature prominently in this paper. In Section~\ref{sspos} we
associate a simplicial poset $\P(\Gamma)$ to an arbitrary torus
graph $\Gamma$; our main result there (Theorem~\ref{tequi})
establishes an isomorphism between the equivariant cohomology of
$\Gamma$ and the \emph{face ring} of~$\P(\Gamma)$. This theorem
continues the series of results identifying the equivariant
cohomology of a smooth toric variety, a (quasi)toric
manifold~\cite{da-ja91}, and a torus
manifold~\cite[Th.~7.5]{ma-pa06} with the face ring of the
appropriate polytope, simplicial complex, or simplicial poset.

Despite the concepts of GKM- and torus graphs diverge in general,
they have an important subclass of \emph{$n$-independent}
GKM-graphs in their intersection. Therefore, our methods and
results about torus graphs are fully applicable to this subclass
of GKM-graphs, which may be considered as a partial answer to some
questions about GKM-graphs posed in the introduction
of~\cite{gu-za01}.

Apart from topological applications to the study of torus action,
the concept of a torus graph and the associated simplicial poset
appears to be of considerable interest for combinatorial
commutative algebra. Since the appearance of Stanley's
book~\cite{stan96} the face ring (or the \emph{Stanley--Reisner
ring}) of simplicial complex has become one of the most important
media of applications of commutative-algebraic methods to
combinatorics. The notion of a face ring has been later extended
to simplicial posets in~\cite{stan91} (see also~\cite{ma-pa06}).
Our primary combinatorial application is the proof of equivalence
of the Cohen--Macaulay properties for simplicial posets and their
face rings. A poset is said to be \emph{Cohen--Macaulay} if the
face ring of its order complex is a Cohen--Macaulay ring. However,
in the case of a simplicial poset $\mathcal P$ the face ring is
defined for the poset $\P$ itself, not only for its order complex.
Therefore, a natural question of whether the Cohen--Macaulay
property can be read directly from the face ring of $\P$ arises.
In~\cite{stan91} Stanley proved that the face ring of a
Cohen--Macaulay simplicial poset is Cohen--Macaulay. In
Theorem~\ref{theo:cmpos} we prove the converse. To do that one has
to show that if the face ring of $\P$ is Cohen--Macaulay, then the
face ring of the order complex of $\P$ is also Cohen--Macaulay.
The passage to the order complex is known to topologists as the
\emph{barycentric subdivision}, and our proof proceeds inductively
by decomposing the barycentric subdivision into a sequence of
elementary \emph{stellar subdivisions} and showing that the
Cohen--Macaulay property is preserved at each step. Stellar
subdivisions of simplicial posets are related to \emph{blow-ups}
of torus manifolds and torus graphs; we further explore this link
in Section~\ref{blowu} by studying the behaviour of equivariant
cohomology under these operations.

In section \ref{spseu} we give a partial answer to the question of
characterising simplicial posets arising from torus graphs. We
also discuss related notions of \emph{orientation} and
\emph{orientability} of a torus graph. Here lies yet another
distinction between the GKM- and torus graphs; all GKM-graphs are
orientable by their definition.

In the last section we deduce certain combinatorial identities for
the number of faces of simplicial posets and torus graphs, which
may be regarded as a yet another generalisation to the
\emph{Dehn--Sommerville equations} for simple polytopes, sphere
triangulations, Eulerian posets etc.

The authors thank the referee for a number of comments and
suggestions, which helped in improving the text significantly.

\section{Torus manifolds and equivariant cohomology}

A \emph{torus manifold} \cite{ha-ma03} is a $2n$-dimensional
compact smooth manifold $M$ with an effective (or faithful) action
of an $n$-dimensional torus $T$ whose fixed point set is
non-empty. This fixed point set $M^T$ is easily seen to consist of
finite number of isolated points. A \emph{characteristic
submanifold} of~$M$ is a codimension-two connected component of
the set fixed pointwise by a circle subgroup of $T$. An
\emph{omniorientation}~\cite{bu-ra01} of $M$ consists of a choice
of orientation for $M$ and for each characteristic submanifold. A
choice of omniorientation allows us to regard the normal bundles
to characteristic submanifolds as complex line bundles, and is
particularly useful for studying the equivariant cohomology,
characteristic classes, stably almost complex structures on torus
manifolds etc. Sometimes fixing an omniorientation is required in
the definition of a torus manifold.

All the cohomology in this paper is taken with $\Z$ coefficients,
unless otherwise specified.

Let $ET\to BT$ be the universal $T$-bundle, with $T$ acting on
$ET$ freely from the right. Let $ET\times_T M$ be the orbit space
of the $T$-action on $ET\times M$ defined by $(u,x)\to
(ug^{-1},gx)$ for $(u,x)\in ET\times M$ and $g\in T$. The
projection onto the first factor gives rise to a fibration
\begin{equation} \label{eqn:fibr}
  M\longrightarrow ET\times_T M \longrightarrow BT
\end{equation}
and the equivariant cohomology of $M$ is defined as the ordinary cohomology
of the total space of this fibration:
\[
  H^*_T(M):=H^*(ET\times_T M).
\]

Assume that $H^{odd}(M)=0$. According to Lemma~2.1
of~\cite{ma-pa06}, this is equivalent to $H^*_T(M)$ being
isomorphic to $H^*(M)\otimes H^*(BT)$ as an $H^*(BT)$-module.
Therefore, $H^*_T(M)$ is a free $H^*(BT)$-module, and the Serre
spectral sequence of the fibration (\ref{eqn:fibr}) collapses.
(This condition is referred to as the \emph{equivariant formality}
of $M$ in~\cite[\S1.1]{gu-za99}, although it is different from the
notion of formality, either plain or equivariant, adopted in the
rational homotopy theory.) Under such an assumption, the
localisation theorem~\cite{hsia75} implies that the restriction
homomorphism
\begin{equation} \label{eqn:rest}
  i^*\colon H^*_T(M)\to H^*_T(M^T)
\end{equation}
is injective, where $i\colon M^T\to M$ is the inclusion. The image
of $H^*_T(M)$ in $H^*_T(M^T)$ can be identified in the same way as
it was done by Goresky, Kottwitz and MacPherson~\cite{g-k-m98} for
their class of manifolds (now known as the \emph{GKM-manifolds}).
We briefly describe their result here. Let $\Sigma_M$ denote the
set of 2-dimensional submanifolds of $M$ each of which is fixed
pointwise by a codimension one subtorus of~$T$. Then every
$S\in\Sigma_M$ is diffeomorphic to a sphere, contains exactly two
$T$-fixed points, and is a connected component of the intersection
of some $n-1$ characteristic submanifolds. Denote by $T_S$ the
isotropy subgroup of~$S$. We have a canonical identification
\[
  H^*_T(M^T)=\Map(M^T,H^*(BT)),
\]
and for each $p\in M^T$ there are exactly $n$ spheres in
$\Sigma_M$ containing~$p$.

\begin{theo}[{\cite{g-k-m98}, \cite[\S11.8]{gu-st99}}, {\cite[Th.~3.1]{h-h-h05}}]\label{theo:gkm}
Suppose $H^*_T(M)$ is a free $H^*(BT)$-module. Then $f\in
\Map(M^T,H^*(BT))$ belongs to the image of the map $i^*$
from~\eqref{eqn:rest} if and only if
\[
  r_S(f(p))=r_S(f(q)) \quad\text{for every 2-sphere $S$ in $\Sigma_M$}
\]
where $r_S$ denotes the restriction $H^*(BT)\to H^*(BT_S)$ and
$p,q$ are the $T$-fixed points in $S$.
\end{theo}
In~\cite{gu-st99} this result is stated for GKM-manifolds and with
coefficients in a field of zero characteristic, but it also holds
for torus manifolds with integer coefficients as stated above
(see~\cite{ma-pa06} and Example~\ref{tmtgi} below). The proof
of~\cite{gu-st99} relies on a result of~Chang and
Skjelbred~\cite{ch-sk74} and the localisation
theorem~\cite{hsia75}. In~\cite{h-h-h05} the theorem is proved
with integer coefficients in a much more general context of
$G$-equivariant cohomology theories, under some additional
assumptions.

The tangential $T$-representation $\tau_pM$ at $p\in M^T$
decomposes into a direct sum of irreducible real two-dimensional
$T$-representations. An omniorientation on $M$ determines
orientations on the corresponding two-dimensional
$T$-representation spaces, so that we may think of them as complex
one-dimensional $T$-representations. Therefore, we have
\begin{equation} \label{eqn:deco}
  \tau_pM=\bigoplus_{i=1}^n V(w_{p,i})
\end{equation}
where $V(w_{p,i})$ denotes a complex one-dimensional
$T$-representations with weight $w_{p,i}$. The set of complex
one-dimensional $T$-representations bijectively corresponds to
$H^2(BT)$. Through this bijection, we may think of an element of
$H^2(BT)$ as a weight of the corresponding $T$-representation.

\section{Torus graphs}\label{stogr}
In their study of GKM-manifolds, Guillemin and Zara~\cite{gu-za99}
introduced a combinatorial object called a \emph{GKM-graph} and
defined a notion of (equivariant) cohomology for such graphs
accordingly. In this section we shall see that a similar idea
works for torus manifolds with a little modification.

Let $\G$ be a connected regular $n$-valent graph, $V(\G)$ the set
of vertices of $\G$, and $E(\G)$ the set of oriented edges of $\G$
(so that each edge of $\G$ enters $E(\G)$ with two possible
orientations). We denote by $i(e)$ and $t(e)$ the initial and
terminal points of $e\in E(\G)$ respectively, and by $\bar e$ the
edge $e$ with the orientation reversed. For $p\in V(\G)$ we set
\[
  E(\G)_p:=\{ e\in E(\G)\mid i(e)=p\}.
\]
A collection $\theta=\{\theta_e\}$ of bijections
\[
  \theta_e \colon E(\Gamma)_{i(e)}\to E(\Gamma)_{t(e)},\qquad e\in E(\G),
\]
is called a \emph{connection} on $\mG$ if
\begin{itemize}
\item[(a)] $\theta_{\bar e}$ is the inverse of $\theta_e$;

\item[(b)] $\theta_e(e)=\bar e$.
\end{itemize}

An $n$-valent graph $\G$ admits $((n-1)!)^{g}$ different
connections, where $g$ is the number of (non-oriented) edges in
$\G$. Slightly modifying the original definition of Guillemin and
Zara~\cite{gu-za99}, we call a map
\[
  \alpha\colon E(\G)\to \Hom(T,S^1)=\Htwo
\]
an \emph{axial function} (associated with the connection $\theta$)
if it satisfies the following three conditions:
\begin{itemize}
\item[(a)] $\alpha(\bar e)=\pm \alpha(e)$;

\item[(b)] elements of $\alpha(E(\Gamma)_p)$ are pairwise linearly
independent (2-\-in\-de\-pen\-dent) for each $p\in V(\G)$;

\item[(c)] $\alpha(\theta_e(e'))\equiv \alpha(e') \mod \alpha(e)$
for any $e\in E(\G)$ and $e'\in E(\Gamma)_{i(e)}$.
\end{itemize}

We also denote $T_e:=\ker\alpha(e)$, the codimension-one subtorus
in $T$ determined by $\alpha$ and~$e$. Then we may reformulate the
condition~(c) above by requiring that the restrictions of
$\alpha(\theta_e(e'))$ and $\alpha(e')$ to $H^*(BT_e)$ coincide.

\begin{rema} Guillemin and Zara required $\alpha(\bar
e)=-\alpha(e)$ in~(a) above. The connection $\theta$ which
satisfies condition~(c) above is unique if elements of
$\alpha(E(\Gamma)_p)$ are 3-independent~\cite{gu-za99}.
\end{rema}

\begin{defi}
We call $\alpha$ a \emph{torus axial function} if it is
$n$-independent, that is, $\alpha(E(\Gamma)_p)$ is a basis of
$H^2(BT)$ for each $p\in V(\G)$. The triple $(\G,\alpha,\theta)$
is then called a \emph{torus graph}. Since the connection $\theta$
is uniquely determined by the second remark above, we may suppress
it in the notation. In what follows we only consider torus axial
functions.
\end{defi}

\begin{rema}
In comparison with the GKM-graphs, the definition of torus graphs
weakens the assumption~(a) on an axial function (by only requiring
$\alpha(\bar e)=\pm\alpha(e)$ instead of $\alpha(\bar
e)=-\alpha(e)$), but strengthens~(b) (by requiring $\alpha$ to be
$n$-independent instead of just 2-independent). Although the
$n$-independence assumption is usually too strict for GKM graphs
(and leaves out some important examples), weakening the other
assumption balances this, as is shown in our next examples.
\end{rema}

\begin{exam}\label{etmtg}
Let $M$ be a torus manifold. Define a regular $n$-valent graph
$\Gamma_M$ whose vertex set is $M^T$ and whose edges correspond to
2-spheres from $\Sigma_M$. The summands in~\eqref{eqn:deco}
correspond to the oriented edges of $\G_M$ having $p$ as the
initial point. We assign $w_{p,i}$ to the oriented edge
corresponding to $V(w_{p,i})$. This gives a function
\[
  \alpha_M\colon E(\G_M)\to \Htwo.
\]
The normal bundle of the 2-sphere corresponding to an oriented
edge in $E(\G)$ decomposes into a Whitney sum of complex $T$-line
bundles. This decomposition defines a connection $\theta_M$ in
$\G_M$. It is not difficult to see that $\alpha_M$ satisfies the
three conditions from the definition of torus axial function.
\end{exam}

\begin{exam}\label{ebana}
Two simple examples of torus graphs $\Gamma$ are shown on
Fig.~\ref{ftogr}. The first is 2-valent and the second is
3-valent. The axial function $\alpha$ assigns the basis elements
$t_1,t_2\in H^2(BT^2)$ (resp. $t_1,t_2,t_3\in H^2(BT^3)$) to the
two (resp. three) edges of $\Gamma$, regardless of the
orientation. These torus graphs are not GKM-graphs, as the
condition $\alpha(\bar e)=-\alpha(e)$ is not satisfied. Both come
from torus manifolds, $S^4$ and $S^6$ respectively, where the
torus action is obtained by suspending the standard coordinatewise
torus actions on $S^3$ and $S^5$ (see~\cite[Ex.~3.2]{ma-pa06}).
\begin{figure}[h]
\begin{picture}(120,30)
\put(30,0){\circle*{1}}
\put(30,30){\circle*{1}}
\put(90,0){\circle*{1}}
\put(90,30){\circle*{1}}
\qbezier(30,0)(15,15)(30,30)
\qbezier(30,0)(45,15)(30,30)
\qbezier(90,0)(70,15)(90,30)
\qbezier(90,0)(110,15)(90,30)
\qbezier(90,0)(95,15)(90,30)
\put(19,14){$t_1$}
\put(38.5,14){$t_2$}
\put(76.5,14){$t_1$}
\put(89,14){$t_2$}
\put(101,14){$t_3$}
\put(22,-5){(a) \ $n=2$}
\put(82,-5){(b) \ $n=3$}
\end{picture}
\medskip
\caption{Torus graphs.}
\label{ftogr}
\end{figure}
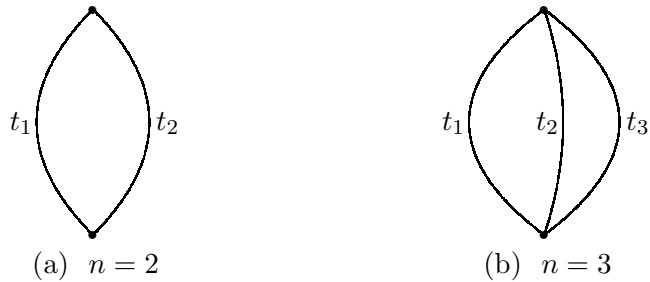
\end{exam}

\begin{defi}
The \emph{equivariant cohomology} $\HT(\mG)$ of a torus graph
$\mG$ is a set of maps
\[
  f\colon V(\G)\to \SV
\]
such that for every $e\in E(\Gamma)$ the restrictions of $f(i(e))$
and $f(t(e))$ to $H^*(BT_e)$ coincide. Since $H^*(BT)$ is a ring,
the vertex-wise multiplication endows the function space
$H^*(BT)^{V(\Gamma)}$ with a ring structure. Its subspace
$H_T^*(\Gamma)$ also becomes a ring because the restriction map
$H^*(BT)\to H^*(BT_e)$ is multiplicative. Moreover,
$H_T^*(\Gamma)$ is an algebra over $\SV$.
\end{defi}

\begin{exam}\label{tmtgi}\sloppy
If $H^*_T(M)$ is a free $H^*(BT)$-module (which happens if
$H^{odd}(M)=0$) and $\mG_M:=(\G_M,\alpha_M,\theta_M)$ is the
associated torus graph, then there is a ring isomorphism
$H^*_T(\mG_M)\cong H^*_T(M)$. This follows either from
Theorem~\ref{theo:gkm} or from the explicit calculation of the two
rings, given in Theorem~\ref{tequi} below
and~\cite[Cor.~7.6]{ma-pa06}.
\end{exam}

\section{Faces and Thom classes}

\begin{defi}
Let $(\G,\alpha,\theta)$ be a torus graph, and $\G'$ a connected
regular $k$-valent subgraph of~$\G$, where $0\le k\le n$. We say
that $(\G',\alpha|E(\G'))$ is a \emph{$k$-dimensional face} of
$\mG$ if $\G'$ is invariant under the connection $\theta$. We
refer to $(n-1)$-dimensional faces as \emph{facets}.
\end{defi}

An intersection of faces is a union of faces. We define the
\emph{Thom class} of a $k$-face $\mF=(\G',\alpha|E(\G'))$ as a map
$\ta{\mF}\colon V(\Gamma)\to H^{2(n-k)}(BT)$ where
\begin{equation}\label{etaf1}
  \ta{F}(p):=\begin{cases}\displaystyle{\prod_{i(e)=p,\ e\notin\Gamma'}
  \alpha(e)}\quad&\text{if $p\in V(\Gamma')$,}\\
  \qquad 0 \quad&\text{otherwise.}
  \end{cases}
\end{equation}

\begin{lemm}
$\ta{F}$ is an element of $\HT(\Gamma)$.
\end{lemm}
\begin{proof}
Let $e\in E(\Gamma)$. If neither vertex of $e$ is contained in
$\mF$, then the values of $\ta{\mF}$ on both vertices of $e$ are
zero. If only one vertex of $e$, say $i(e)$, is contained in
$\mF$, then $\ta{\mF}(t(e))=0$, while
$\ta{\mF}(i(e))=0\mod\alpha(e)$, so that the restriction of
$\ta{\mF}(i(e))$ to $H^*(BT_e)$ is also zero. Finally, assume that
the whole $e$ is contained in $\mF$. Let $e'$ be an edge such that
$i(e')=i(e)$ and $e'\notin\mF$, so that $\alpha(e')$ is a factor
in $\ta{\mF}(i(e))$. Since $\mF$ is invariant under the
connection, we have $\theta_{e}(e')\notin\mF$. Therefore,
$\alpha(\theta_e(e'))$ is one of the factors in $\ta{\mF}(t(e))$.
Now we have $\alpha(\theta_e(e'))\equiv\alpha(e')\mod\alpha(e)$ by
the definition of axial function. The same relation holds for
every other factor in $\ta{\mF}(i(e))$, whence the restrictions of
$\ta{\mF}(i(e))$ and $\ta{\mF}(t(e))$ to $H^*(BT_e)$ coincide.
\end{proof}

\begin{lemm} \label{lemm:face}
If $\mG$ is a torus graph, then there is a unique $k$-face of
$\mG$ containing any given $k$ elements in $E(\G)_p$.
\end{lemm}
\begin{proof}
Let $W$ be the $k$-dimensional subspace of $H^2(BT)$ spanned by
the images of the given $k$ edges in $E(\G)_p$ under the axial
function~$\alpha$. Take any element $e$ from $E(\G)_p$.  Through
the connection $\theta_e$ the given $k$ edges in $E(\G)_p$ map to
some $k$ edges in $E(\G)_{t(e)}$. The $\alpha$-images of these $k$
edges in $E(\G)_{t(e)}$ span the same subspace $W$ in $H^2(BT)$.
Proceeding in the same way, we translate the given $k$ edges in
$E(\G)_p$ along the edges related to $E(\G)_p$ via the connection.
The uniqueness of the connection guarantees that the resulting
graph is regular and $k$-valent.
\end{proof}

The intersection of two faces $G$ and $H$ of $\Gamma$ is a finite
set of faces. Denote by $G\vee H$ a minimal face containing both
$G$ and $H$. In general such a least upper bound may fail to exist
or be non-unique; however it exists and is unique provided that
the intersection $G\cap H$ is non-empty.

\begin{lemm}\label{ltcre}
For any two faces $G$ and $H$ of $\Gamma$ the corresponding Thom
classes satisfy the relation
\begin{equation}\label{etocl}
  \ta{G}\ta{H}=
    \ta{G\vee H}\cdot\!\!\!\sum_{E\in{G\cap H}}\!\!\!\ta{E},
\end{equation}
where we formally set $\ta{\Gamma}=1$ and $\tau_{\varnothing}=0$,
and the sum in the right hand side is taken over connected
components $E$ of $G\cap H$.
\end{lemm}
\begin{proof} Although the proof is the same as that
of~\cite[Lemma~6.3]{ma-pa06}, we include it here for reader's
convenience. Take $p\in V(\Gamma)$. For a face $F$ such that $p\in
F$, we set
\[
  N_p(F):= \{ e\in E(\Gamma)\colon i(e)=p,\ e\notin F\},
\]
which may be thought of as the set of directions normal to $F$
at~$p$. Then the identity \eqref{etaf1} can be written as
\begin{equation} \label{etaf2}
  \ta{F}(p)=\prod_{e\in N_p(F)}\alpha(e)
\end{equation}
where the right hand side is understood to be 1 if
$N_p(F)=\varnothing$ and $0$ if $p\notin F$. If $p\notin G\cap H$,
then $p\notin E$ for any connected component $E$ of $G\cap H$ and
either $p\notin G$ or $p\notin H$. Therefore, both sides
of~\eqref{etocl} take zero value on~$p$. If $p\in G\cap H$, then
\[
  N_p(G)\cup N_p(H)=N_p(G\vee H)\cup N_p(E)
\]
where $E$ is the connected component of $G\cap H$ containing $p$,
and $p\notin E'$ for any other connected component $E'\in G\cap
H$. This together with \eqref{etaf2} shows that both sides
of~\eqref{etocl} take the same value on~$p$.
\end{proof}

\begin{lemm}\label{ltcmg}
The Thom classes $\ta{F}$ corresponding to all proper faces of
$\Gamma$ constitute a set of ring generators for $H^*_T(\Gamma)$.
\end{lemm}
\begin{proof}
Again, the proof is very similar to that
of~\cite[Prop.~7.4]{ma-pa06}. Let $\eta\in H^{>0}_T(\Gamma)$ be a
nonzero element. Set
\[
  Z(\eta):=\{p\in V(\Gamma)\colon \eta(p)=0\}.
\]
Take $p\in V(\Gamma)$ such that $p\notin Z(\eta)$. Then
$\eta(p)\in H^*(BT)$ is non-zero and we can express it as a
polynomial in $\{\alpha(e)\colon e\in E(\Gamma)_p\}$, which is a
basis of $H^2(BT)$. Let
\begin{equation}
\label{monom}
  \prod_{e\in E(\Gamma)_p}\alpha(e)^{n_e},\qquad n_e\ge0,
\end{equation}
be a monomial entering $\eta(p)$ with a non-zero coefficient. Let
$F$ be the face spanned by the edges $e$ with $n_e=0$. Denote by
$I(F)$ the ideal in $H^*(BT)$ generated by all elements
$\alpha(e)$ with $e\in F$. Then $\eta(p)\notin I(F)$ since
$\eta(p)$ contains monomial~\eqref{monom}. Suppose $\eta(q)\in
I(F)$ for some other vertex $q\in F$. Then $\eta(s)\in I(F)$ for
any vertex $s\in F$ joined to $q$ by an edge $f\subseteq F$
because $\eta(q)-\eta(s)$ is divisible by $\alpha(f)$ by the
definition of axial function. Since $F$ is a connected subgraph,
$\eta(q)\in I(F)$ for any vertex $q\in F$, in contradiction with
$\eta(p)\notin I(F)$. Hence, $\eta(q)\notin I(F)$, in particular
$\eta(q)\ne0$, for every vertex $q\in F$.

On the other hand, it follows from~\eqref{etaf1} that
monomial~\eqref{monom} can be written as $\u{F}\ta{F}(p)$ where
$\u{F}$ is a product of some Thom classes corresponding to faces
containing~$F$. Set $\eta':=\eta-\u{F}\ta{F}\in H^*_T(\Gamma)$.
Since $\ta{F}(q)=0$ for $q\notin F$, we have $\eta'(q)=\eta(q)$
for all such~$q$. At the same time, $\eta(q)\ne0$ for every vertex
$q\in F$ by the argument from the previous paragraph. It follows
that $Z(\eta')\supseteq Z(\eta)$. However, the number of monomials
in $\eta'(p)$ is less than that in $\eta(p)$. Therefore,
subtracting from $\eta$ monomials in Thom classes we can
eventually achieve an element $\lambda$ such that $Z(\lambda)$
contains $Z(\eta)$ as a proper subset. Repeating this procedure,
we end up at an element whose value on every vertex is zero.
\end{proof}

In order to finish our description of the equivariant cohomology
of torus graphs we need a combinatorial diversion to the concepts
of simplicial posets and face rings.

\section{Simplicial posets}\label{sspos}
We start by briefly reviewing simplicial posets and related
algebraic notions. Then we prove our main result here,
Theorem~\ref{tequi}, which effectively describes the equivariant
cohomology of torus graphs. The discussion of simplicial posets
continues in the next section, where we concentrate on the
Cohen--Macaulay property.

A poset $\mathcal P$ is called \emph{simplicial} if it has an
initial element $\hatzero$ and for each $\sigma\in\mathcal P$ the
lower segment $[\hatzero,\sigma]$ is a boolean lattice (the face
poset of a simplex). We assume all our posets to be finite. An
example of a simplicial poset is provided by the face poset of a
simplicial complex, but there are many simplicial posets that do
not arise in this way. We identify a simplicial complex with its
face poset, thereby regarding simplicial complexes as particular
cases of simplicial posets.

To each $\sigma\in\mathcal P$ we assign a geometric simplex whose
face poset is $[\hatzero,\sigma]$, and glue these geometrical
simplices together according to the order relation in $\mathcal
P$. We get a cell complex in which the closure of each cell can be
identified with a simplex preserving the face structure and all
the attaching maps are inclusions. We call it a \emph{simplicial
cell complex} and denote its underlying space by $|\mathcal P|$.
In what follows we shall not distinguish between simplicial posets
and simplicial cell complexes and refer to elements
$\sigma\in\mathcal P$ as \emph{simplices}.

Let $\mathcal P$ be a simplicial poset. The \emph{rank function}
on $\P$ is defined by setting $\mathop{\mathrm{rk}}\sigma=k$ for
$\sigma\in\mathcal P$ if $[\hatzero,\sigma]$ is the face poset of
a $(k-1)$-dimensional simplex. The rank of $\mathcal P$ is the
maximum of ranks of its elements. Let $\k$ be a commutative ring
with unit. Introduce the graded polynomial ring $\k[v_\sigma\colon
\sigma\in\mathcal P\setminus\hatzero]$ with $\deg
v_\sigma=2\mathop{\mathrm{rk}}\sigma$. We also write formally
$v_{\hatzero}=1$. For any two simplices $\sigma,\tau\in\mathcal P$
denote by $\sigma\vee\tau$ the set of their least common upper
bounds (\emph{joins}), and by $\sigma\wedge\tau$ the set of their
greatest common lower bounds (\emph{meets}). Since $\mathcal P$ is
a simplicial poset, $\sigma\wedge\tau$ consists of a single
simplex whenever $\sigma\vee\tau$ is non-empty. There is the
following simple characterisation of the subclass of simplicial
complexes in simplicial posets.

\begin{prop}\label{prop:psico}
$\P$ is a simplicial complex if and only if for any two elements
$\sigma,\tau\in\P$ the set $\sigma\vee\tau$ is either empty or
consists of a single simplex.
\end{prop}
\begin{proof}
Let $V(\P)=\{v_1,\ldots,v_m\}$ be the set of vertices (rank one
elements) of~$\P$. Introduce a simplicial complex $K$ on the
vertex set $V(\P)$ whose simplices are those subsets
$\{v_{i_1},\ldots,v_{i_k}\}$ for which there is an element
$\sigma\in\P$ with such vertex set. There is an obvious surjective
order preserving map $\P\to K$ assigning to an element of $\P$ its
vertex set. The injectivity of this map follows from the
additional assumption on the joins. Therefore, $\P$ is (the face
poset of)~$K$. The other direction is obvious.
\end{proof}

\begin{defi}[\cite{stan91}]
The \emph{face ring} of a simplicial poset $\mathcal P$ is the
quotient
\[
  \k[\mathcal P]:=
  \k[v_\sigma\colon \sigma\in\mathcal P\setminus\hatzero]/\mathcal I_{\mathcal P},
\]
where $\mathcal I_{\mathcal P}$ is the ideal generated by the
elements
\begin{equation}\label{estre}
  v_\sigma v_\tau-v_{\sigma\wedge\tau}\cdot\!\!\!\sum_{\eta\in\sigma\vee\tau}\!\!\!v_\eta.
\end{equation}
The sum over the empty set is assumed to be zero, so we have
$v_\sigma v_\tau=0$ if $\sigma\vee\tau=\varnothing$.
\end{defi}

\begin{rema}
The definition above extends the notion of the face ring of a
simplicial complex (also known as the \emph{Stanley--Reisner
ring}) to simplicial posets. In the case when $\mathcal P$ is a
simplicial complex we may rewrite~\eqref{estre} as $v_\sigma
v_\tau-v_{\sigma\wedge\tau}v_{\sigma\vee\tau}$ (because
$\sigma\vee\tau$ is either empty or consists of a single simplex),
and use the latter relation to express every element $\sigma\in\P$
as
\[
  v_\sigma=\prod_{v_i\in V(\sigma)}v_i,
\]
where $V(\sigma)$ is the vertex set of $\sigma$. The relations
between $v_i$'s coming from~\eqref{estre} can now be written as
\begin{equation}\label{screl}
v_{i_1}\cdots v_{i_k}=0\quad\text{if
}\{v_{i_1},\ldots,v_{i_k}\}\text{ does not span a simplex of }\P.
\end{equation}
The face ring $\k[\P]$ is isomorphic to the quotient of the
polynomial ring $\k[v_1,\ldots,v_m]$ by~\eqref{screl}, where
$V(\P)=\{v_1,\ldots,v_m\}$ and $\deg v_i=2$. This is the standard
form of the face ring of a simplicial complex~\cite{stan96}.
\end{rema}

We briefly remind several algebraic constructions
from~\cite{ma-pa06}.

\begin{lemm}[{\cite[Lemma~5.4]{ma-pa06}}]\label{chain}
Every element of $\k[\P]$ can be uniquely written as a linear
combination of monomials
$v_{\tau_1}^{\alpha_1}v_{\tau_2}^{\alpha_2}\cdots
v_{\tau_k}^{\alpha_k}$ corresponding to chains of fully ordered
elements $\tau_1<\tau_2<\ldots<\tau_k$ of~$\P$.
\end{lemm}

In other words, the monomials
$v_{\tau_1}^{\alpha_1}v_{\tau_2}^{\alpha_2}\cdots
v_{\tau_k}^{\alpha_k}$ with $\tau_1<\tau_2<\ldots<\tau_k$
constitute an additive basis of $\k[\P]$. We refer to the
expansion of an element $x\in\k[\P]$ in terms of this basis as the
\emph{chain decomposition} of~$x$. To achieve a chain
decomposition we inductively use \emph{straightening
relation}~\eqref{estre}, which allows us to express the product of
two unordered elements via the products of ordered elements. This
can be restated by saying that the face ring is an example of an
\emph{algebra with straightening law} (see discussion
in~\cite[p.~323]{stan91}).

Given an element $\sigma\in\P$, define the \emph{restriction map}
$s_\sigma$ as
\[
  s_\sigma\colon\k[\P]\to\k[\P]/(v_\tau\colon\tau\not\le\sigma).
\]
Its codomain is isomorphic to a polynomial ring on $\rank\sigma$
generators. The following is a key algebraic statement, which has
several geometric interpretations.

\begin{lemm}[{\cite[Lemma~5.6]{ma-pa06}}]\label{lemm:algrm}
The sum $s=\bigoplus_{\sigma\in\P}s_\sigma$ of the restrictions to
all elements of~$\P$,
\[
  s\colon\k[\P]\to\bigoplus_{\sigma\in\P}\k[\P]/(v_\tau\colon\tau\not\le\sigma),
\]
is a monomorphism.
\end{lemm}

It is clear that to get a monomorphism it is enough to take the
sum of restrictions to the maximal elements only. The proof of the
above lemma uses the chain decomposition.

Now let $\mG$ be a torus graph. By Lemma~\ref{lemm:face}, the
faces of $\mG$ form a simplicial poset of rank $n$ with respect to
the reversed inclusion relation. We denote this simplicial poset
by $\mathcal P(\mG)$. In section~\ref{blowu} we shall discuss
which simplicial posets may arise in this way. We prefer to stick
to the original face-inclusions notation while dealing with torus
graphs; then the face ring $\k[\mathcal P(\mG)]$ is a quotient of
the polynomial ring on generators $\v{F}$, where $F$ is a proper
non-empty face of $\mG$, and $\deg\v{F}=2(n-\dim F)$. We set
formally $v_\varnothing=0$ and $\v{\mG}=1$; then the defining
relation for the face ring is the same as~\eqref{etocl}.

\begin{exam}
1. Let $\Gamma$ be the torus graph shown on Fig.~\ref{ftogr}~(a),
see Example~\ref{ebana}. Denote its two edges by $e$ and $g$, and
the two vertices by $p$ and $q$. The simplicial cell complex
$\mathcal P(\Gamma)$ is obtained by gluing two segments along
their boundaries. (It looks the same as $\Gamma$ itself, but this
is a mere coincidence, see the second example below.) The face
ring $\k[\mathcal P(\mG)]$ is the quotient of the graded
polynomial ring
\[
  \k[v_e,v_g,v_p,v_q],\quad\deg v_e= \deg v_g=2,\quad
  \deg v_p=\deg v_q=4
\]
by the two relations
\[
  v_ev_g=v_p+v_q,\quad v_pv_q=0.
\]

2. Now let $\Gamma$ be the torus graph shown on
Fig.~\ref{ftogr}~(b). Denote its vertices by $p$ and $q$, the
edges by $e$, $g$, $h$, and the 2-faces by $E$, $G$, $H$ so that
$e$ is opposite to $E$, etc. The simplicial cell complex $\mathcal
P(\Gamma)$ is obtained by gluing two triangles along their
boundaries. The face ring $\k[\mathcal P(\mG)]$ is isomorphic to
the quotient of the graded polynomial ring
\[
  \k[\v{E},\v{G},\v{H},v_p,v_q],\quad\deg\v{E}=\deg\v{G}=\deg\v{H}=2,\quad
  \deg v_p=\deg v_q=6
\]
by the two relations
\[
  \v{E}\v{G}\v{H}=v_p+v_q,\quad v_pv_q=0.
\]
(The generators corresponding to the edges can be excluded because
of the relations $v_e=\v{G}\v{H}$, etc.)
\end{exam}

By the definition of the equivariant cohomology of torus graph, it
comes together with a monomorphism into the sum of polynomial
rings:
\[
  r\colon\HT(\mG)\longrightarrow\bigoplus_{V(\Gamma)}H^*(BT),
\]
whose analogy with the algebraic restriction map $s$ from
Lemma~\ref{lemm:algrm} now becomes clear. The latter can now be
written as
\[
  s\colon\k[\mathcal P(\mG)]\longrightarrow\bigoplus_{p\in V(\Gamma)}
  \k[\mathcal P(\mG)]/(\v{F}\colon F\not\ni p).
\]

\begin{theo}\label{tequi}
$H^*_T(\Gamma)$ is isomorphic to the face ring $\Z[\mathcal
P(\mG)]$. In other words, $H^*_T(\Gamma)$ is isomorphic to the
quotient of the graded polynomial ring generated by the Thom
classes $\ta{F}$ modulo relations~\eqref{etocl}.
\end{theo}
\begin{proof}
We start by constructing a map
\[
  \Z[\v{F}\colon F\text{ a face}]\longrightarrow\HT(\mG)
\]
that sends $\v{F}$ to $\ta{F}$. By Lemma~\ref{ltcre}, it factors
through a map $\varphi\colon\Z[\mathcal P(\mG)]\to\HT(\mG)$. This
map is surjective by Lemma~\ref{ltcmg}. Finally, $\varphi$ is
injective, because we have a decomposition $s=r\circ\varphi$, and
$s$ is injective by Lemma~\ref{lemm:algrm}.
\end{proof}

\section{Cohen--Macaulay rings, complexes, and posets}
A simplicial complex $K$ is called \emph{Cohen--Macaulay} (over
$\k$) if its face ring $\k[K]$ is a Cohen--Macaulay ring (see,
e.g.,~\cite{stan96}). The Cohen--Macaulay property has several
topological and algebraic interpretations (some of which we list
below) and is very important for both topological and
combinatorial applications of the face rings.

We shall not give a definition of a Cohen--Macaulay ring in the
general case; instead we state a proposition characterising
Cohen--Macaulay face rings of simplicial complexes.

\begin{prop}[{see~\cite[Ch.~II]{stan96}
or~\cite[Ch.~3]{bu-pa02}}]\label{prop:rlsop} Assume $K$ is an
$(n-1)$-dimensional simplicial complex and $\k$ is an infinite
field. Then $\k[K]$ is a Cohen--Macaulay ring if and only if there
exists a sequence $\theta_1,\ldots,\theta_n$ of linear (i.e.,
degree-two) elements of $\k[K]$ satisfying one of the two
following equivalent conditions:
\begin{itemize}
\item[(a)] $\theta_i$ is not a zero divisor in
$\k[K]/(\theta_1,\ldots,\theta_{i-1})$ for $i=1,\ldots,n$;

\item[(b)] $\theta_1,\ldots,\theta_n$ are algebraically
independent and $\k[K]$ is a free finitely generated module over
its polynomial subring $\k[\theta_1,\ldots,\theta_n]$.
\end{itemize}
\end{prop}

A sequence satisfying the first condition above is called
\emph{regular}. A sequence $\theta_1,\ldots,\theta_n$ of
algebraically independent linear elements for which $\k[K]$ is a
finitely generated module over $\k[\theta_1,\ldots,\theta_n]$
(i.e., $\k[K]/(\theta_1,\ldots,\theta_n)$ is a finite dimensional
$\k$-vector space) is called an \emph{lsop} (linear systems of
parameters). Remember that $\dim K=n-1$, so that an lsop in
$\k[K]$ must have length~$n$. An lsop always exists over an
infinite field (see~\cite[Th.~1.5.17]{br-he98}); the existence of
an lsop over a finite field or $\Z$ is usually a subtle issue.
Thus, we may reformulate Proposition~\ref{prop:rlsop} by saying
that $K$ is Cohen--Macaulay (over an infinite $\k$) if and only if
$\k[K]$ admits a regular lsop. If $\k[K]$ is Cohen--Macaulay, then
every lsop is regular~\cite[Th.~2.1.2]{br-he98}.

Linear systems of parameters in face rings may be detected by
means of the following result.

\begin{prop}[{see~\cite[Lemma~III.2.4]{stan96}}]\label{prop:dlsop}
A sequence of linear elements $\theta_1,\ldots,\theta_n$ of
$\k[K]$ is an lsop if and only if for every simplex $\sigma\in K$
the images $s_\sigma(\theta_1),\ldots,s_\sigma(\theta_n)$ under
the restriction map
\[
  s_\sigma\colon\k[K]\to\k[K]/(v_i\colon i\notin\sigma)
\]
generate the positive degree part of the polynomial ring
$\k[K]/(v_i\colon i\notin\sigma)$.
\end{prop}

A theorem due to Reisner gives a purely topological
characterisation of Cohen--Macaulay complexes. We shall use the
following version of Reisner's theorem, which is due to Munkres.

\begin{theo}[{\cite[Cor.~3.4]{munk84}}]\label{theo:munkr}
A simplicial complex $K$ is Cohen--Macaulay if and only if the
space $X=|K|$ satisfies
\[
  \widetilde{H}^i(X)=0=\widetilde H^i(X,X\setminus p)
\]
for any $p\in X$ and $i<\dim X$.
\end{theo}

Now let $\P$ be a simplicial poset. Its \emph{barycentric
subdivision} $\P'$ is the order complex $\Delta(\overline\P)$ of
the poset $\overline\P=\P\setminus\hatzero$. By the definition,
$\P'$ is a genuine simplicial complex. Its geometric realisation
can be obtained from the simplicial cell complex $\P$ by
barycentrically subdividing each of its simplices.

Following Stanley~\cite{stan91}, we call a simplicial poset $\P$
\emph{Cohen--Macaulay} (over $\k$) if $\P'$ is a Cohen--Macaulay
simplicial complex, that is, if the face ring $\k[\P']$ is
Cohen--Macaulay. A question arises whether the Cohen--Macaulay
property can be read directly from its face ring $\k[\P]$. By a
result of Stanley~\cite[Cor.~3.7]{stan91}, the face ring of a
Cohen--Macaulay simplicial poset is Cohen--Macaulay. The rest of
this section is dedicated to the proof of the converse of this
statement (see Theorem~\ref{theo:cmpos}).

We call a simplicial subdivision of $\P$ \emph{regular} if it is a
genuine simplicial complex. For example, the barycentric
subdivision is a regular subdivision. The following
characterisation of Cohen--Macaulay simplicial posets follows from
Theorem~\ref{theo:munkr}.

\begin{coro}\label{coro:cmcha}
The following conditions are equivalent for $\P$ to be a
Cohen--Macaulay poset:
\begin{itemize}
\item[(a)] the barycentric subdivision of $\P$ is a
Cohen--Macaulay complex;

\item[(b)] any regular subdivision of $\P$ is a Cohen--Macaulay
complex;

\item[(c)] a regular subdivision of $\P$ is a Cohen--Macaulay
complex.
\end{itemize}
\end{coro}

As a further corollary we obtain that Theorem~\ref{theo:munkr}
itself holds for arbitrary simplicial poset, i.e., the
Cohen--Macaulay property for simplicial cell complexes is also of
purely topological nature. All algebraic results from the
beginning of this section also directly generalise to simplicial
posets (for Proposition~\ref{prop:dlsop}
see~\cite[Th.~5.4]{bu-pa04}).

The barycentric subdivision $\P'$ can be obtained as the result of
a sequence of \emph{stellar subdivisions} of~$\P$. Fix a
$(k-1)$-dimensional simplex $\sigma\in\P$. The \emph{star} of
$\sigma$, its \emph{boundary}, and \emph{link} are the following
subposets:
\begin{align*}
  \st_{\P}\sigma&=\{\tau\in\P\colon\sigma\vee\tau\ne\varnothing\},\\
  \partial\st_{\P}\sigma&=\{\tau\in\P\colon\sigma\not\le\tau,\,
  \sigma\vee\tau\ne\varnothing\},\\
  \lk_{\P}\sigma&=\{\tau\in\P\colon\tau\wedge\sigma=\hatzero,\,
  \sigma\vee\tau\ne\varnothing\}.
\end{align*}
These correspond to the usual notions of star (or combinatorial
neighbourhood) of a simplex in a triangulation, its boundary, and
link.

\begin{rema}
The star of a simplex can be thought of as its ``closed
combinatorial neighbourhood''. If $\P$ is (the poset of faces of)
a simplicial complex, then the poset $\lk\sigma$ is isomorphic to
the upper interval
\[
  \P_{>\sigma}=\{\rho\in\P\colon\rho>\sigma\},
\]
and $\st\sigma=\sigma\mathbin{*}\lk\sigma$ (here $*$ denotes the
\emph{join} of simplicial complexes). However this is not the case
in general, see Example~\ref{2tria} below.
\end{rema}

\begin{defi}
Let $\P$ be a simplicial poset and $\sigma\in\P$ a simplex. Assume
first that $\P$ is a simplicial complex. Then the \emph{stellar
subdivision} of $\P$ at $\sigma$ is a simplicial complex
$\widetilde\P$ obtained by removing from $\P$ the star of $\sigma$
and adding the cone over the boundary of the star:
\begin{equation}\label{stsub}
  \widetilde\P=(\P\setminus\st_{\P}\sigma)\cup
  \mathop{\mathrm{cone}}(\partial\st_{\P}\sigma).
\end{equation}
Therefore, if $v$ is the new vertex of~$\widetilde\P$, then we
have $\lk_{\widetilde\P}v=\partial\st_{\P}\sigma$ and
$|\st_{\P}\sigma|\cong|\st_{\widetilde\P}v|$. Now, if $\P$ is an
arbitrary simplicial poset, then its stellar subdivision
$\widetilde\P$ at $\sigma$ is obtained by stellarly subdividing
each simplex containing~$\sigma$. The term ``subdvision'' is
justified by the fact that $\P$ and $\widetilde\P$ are
homeomorphic as topological spaces.
\end{defi}

\begin{prop}\label{prop:barst}
The barycentric subdivision $\P'$ can be obtained as a sequence of
stellar subdivisions of~$\P$. Moreover, each stellar
subdivision in the sequence is taken at a simplex whose star is a
simplicial complex.
\end{prop}
\begin{proof}
Assume $\dim\P=n-1$. We start by taking stellar subdivisions of
all $(n-1)$-dimensional simplices. Denote the resulting complex by
$\P_1$. Then we take stellar subdivisions of $\P_1$ at all
$(n-2)$-dimensional simplices corresponding to $(n-2)$-simplices
of $\P$, and denote the result by $\P_2$. Proceeding this way, at
the end we get a simplicial poset $\P_{n-1}$, which is obtained by
stellar subdivisions of $\P_{n-2}$ at all 1-simplices
corresponding to the edges of~$\P$. Then $\P_{n-1}$ is~$\P'$. To
prove the second statement, assume that $\mathcal R$ and
$\widetilde{\mathcal R}$ are the two subsequent complexes in the
sequence, and $\widetilde{\mathcal R}$ is obtained from $\mathcal
R$ by a stellar subdivision at $\sigma$. Then $\st_{\mathcal
R}\sigma$ is isomorphic to $\sigma\mathbin{*}(\P_{>\sigma})'$ and
thereby is a simplicial complex.
\end{proof}

We proceed with two lemmas necessary to prove our main result.

\begin{lemm}\label{lemm:frppt}
Let $\P$ be a $(n-1)$-dimensional simplicial poset with the vertex
set $V(\P)=\{v_1,\ldots,v_m\}$, and assume that the first $k$
vertices $v_1,\ldots,v_k$ span a face $\sigma$. Assume further
that $\st_\P\sigma$ is a simplicial complex, and consider the
stellar subdivision $\widetilde\P$ of $\P$ at $\sigma$. Let $v$
denote the degree-two generator of $\k[\widetilde\P]$
corresponding to the added vertex. Then there exists a unique map
$\beta\colon\k[\P]\to\k[\widetilde\P]$ such that
\begin{align*}
  v_\tau & \mapsto v_\tau &&\text{if \ }\tau\notin\st_{\P}\sigma,\\
  v_i & \mapsto v+v_i, && i=1,\ldots,k,\\
  v_i & \mapsto v_i, && i=k+1\ldots,m
\end{align*}
(we use the same notation for the vertices and the corresponding
degree-two generators of the face ring). Moreover, $\beta$ is a
monomorphism, and if $\theta_1,\ldots,\theta_n$ is an lsop in
$\k[\P]$, then $\beta(\theta_1),\ldots,\beta(\theta_n)$ is an lsop
in $\k[\widetilde\P]$.
\end{lemm}
\begin{proof}
In order to define the map $\beta$ completely we have to specify
the images of $v_\tau$ for $\tau\in\st_{\P}\sigma$. Choose such a
$v_\tau$ and let $V(\tau)=\{v_{i_1},\ldots,v_{i_\ell}\}$ be its
vertex set. Then we have the following identity in
$\k[\P]=\k[v_\tau\colon\tau\in\P\setminus\hatzero]/\mathcal
I_{\P}$:
\begin{equation}\label{expan}
  v_{i_1}\cdots v_{i_\ell}=
  v_\tau+\sum_{\eta\colon V(\eta)=V(\tau),\,\eta\ne\tau}v_\eta.
\end{equation}
For every $v_\eta$ in the latter sum we have
$\eta\notin\st_{\P}\sigma$ since $\st_{\P}\sigma$ is a simplicial
complex (see Proposition~\ref{prop:psico}). Since $\beta$ is
already defined on the product in the left hand side and on the
sum in the right hand side, this uniquely determines
$\beta(v_\tau)$. Therefore, the map $\beta$ is defined on all
monomials described in Lemma~\ref{chain}, and we may construct a
map of $\k$-modules $\beta\colon\k[\P]\to\k[\widetilde\P]$ using
the chain decomposition.

Next we have to check that $\beta$ is a ring homomorphism.
Consider the projection
\[
  p\colon \k[\P]\to\k[\P]/(v_\tau\colon\tau\notin\st_{\P}\sigma)
  =\k[\st_\P\sigma]
\]
and denote its kernel by $R$. Similarly, denote $\widetilde
R=\ker(\tilde p\colon\k[\widetilde\P]\to\k[\st_{\widetilde\P}v])$.
The ideal $R$ has a $\k$-basis consisting of monomials
$v_{\tau_1}^{\alpha_1}v_{\tau_2}^{\alpha_2}\cdots
v_{\tau_k}^{\alpha_k}$ with $\tau_1<\tau_2<\ldots<\tau_k$, \
$\alpha_i>0$ for $1\le i\le k$ and $\tau_k\notin\st_\P\sigma$.
Since the simplicial cell complexes $\P$ and $\widetilde\P$ do not
differ on the complement to $\st_\P\sigma$ and
$\st_{\widetilde\P}v$ respectively, the map $\beta$ restricts to
the identity isomorphism $R\to\widetilde R$.

The map $\beta$ induces an additive map
$\k[\st_\P\sigma]\to\k[\st_{\widetilde\P}v]$, and our next
observation will be that the latter is a ring homomorphism. Since
$\k[\st_\P\sigma]$ is generated in degree two, we need to check
that $\beta$ vanishes on monomials~\eqref{screl}, that is, that
$\beta(\mathcal I_{\st_\P\sigma})\subset\mathcal
I_{\st_{\widetilde\P}v}$. We may assume that
$\{v_{i_1},\ldots,v_{i_\ell}\}$ is a minimal non-simplex of
$\st_\P\sigma$, that is, every its proper subset is a simplex.
Then we have $\{v_{i_1},\ldots,v_{i_\ell}\}\cap
V(\sigma)=\varnothing$ by the definition of the star. Therefore,
$\beta(v_{i_1}\cdots v_{i_\ell})=v_{i_1}\cdots v_{i_\ell}$, which
belongs to~$\mathcal I_{\st_{\widetilde\P}v}$.

Now we have the following diagram with exact rows:
\begin{equation}\label{exrow}
\begin{CD}
  0 @>>> R @>>> \k[\P] @>p>>
  \k[\st_{\P}\sigma] @>>> 0\\
  @. @VV\cong V @VV\beta V @VVV\\
  0 @>>> \widetilde R @>>> \k[\widetilde\P] @>\tilde p>>
  \k[\st_{\widetilde\P}v] @>>>
  0,
\end{CD}
\end{equation}
in which the left and right vertical arrows are ring
homomorphisms. We need to check that
$\beta(x_1x_2)=\beta(x_1)\beta(x_2)$ for every $x_1,x_2\in\k[\P]$.
Since $\k[\P]=R\oplus\k[\st_{\P}\sigma]$ as $\k$-modules, we may
write $x_i=r_i+s_i$ with $r_i\in R$, \ $s_i\in\k[\st_{\P}\sigma]$,
\ $i=1,2$. For every $s\in\k[\st_{\P}\sigma]$ we have
$\beta(s)=s+vx$ for some $x\in\k[\st_{\widetilde\P}v]$, and $rv=0$
in $\k[\widetilde\P]$ for every $r\in\widetilde R$. Note also that
$rs\in R$ for every $r\in R$, \ $s\in S$, as $R$ is an ideal. Then
we have
\[
  \beta(x_1x_2)=\beta(r_1r_2+r_1s_2+r_2s_1+s_1s_2)=
  r_1r_2+r_1s_2+r_2s_1+\beta(s_1s_2),
\]\\[-15pt]
and\\[-15pt]
\[
  \beta(x_1)\beta(x_2)=(r_1+\beta(s_1))(r_2+\beta(s_2))=
  r_1r_2+r_1\beta(s_2)+r_2\beta(s_1)+\beta(s_1)\beta(s_2).
\]
Since $r_1\beta(s_2)=r_1(s_2+vx_2)=r_1s_2$, \
$r_2\beta(s_1)=r_2s_1$ and $\beta(s_1s_2)=\beta(s_1)\beta(s_2)$,
we conclude that $\beta(x_1x_2)=\beta(x_1)\beta(x_2)$. Thus,
$\beta$ is a ring homomorphism.

The rest of the statement follows by considering the commutative
diagram of restriction maps (see Lemma~\ref{lemm:algrm})
\begin{equation}
\label{epore}
\begin{CD}
  \k[\P] @>\beta>> \k[\widetilde\P]\\
  @VsVV @VsVV\\
  {\displaystyle\bigoplus_{\zeta\in\P}\k[\P]/(v_\tau\colon\tau\not\le\zeta)}
  @>s(\beta)>>
  {\displaystyle\bigoplus_{\zeta\in\widetilde
  \P}\k[\widetilde\P]/(v_\tau\colon\tau\not\le\zeta).}
\end{CD}
\end{equation}
The map $s(\beta)$ sends each direct summand of its domain
isomorphically to at least one summand of its codomain, and
therefore, is a monomorphism (its exact form can be easily guessed
from the definition of $\beta$). Thus, $\beta$ is also a
monomorphism. Finally, the statement about lsop follows from the
above diagram and Proposition~\ref{prop:dlsop}.
\end{proof}

Note that if we map $v_i$ identically for $i=1,\ldots,k$ in the
lemma above, then the map $\k[\P]\to\k[\widetilde\P]$ would still
exist, but fail to be a monomorphism.

\begin{exam}\label{2tria}
The assumption in Lemma~\ref{lemm:frppt} is not satisfied if we
take $\P$ to be the simplicial cell complex obtained by
identifying two 2-simplices along their boundaries and make a
stellar subdivision at a 1-simplex (the star of a 1-simplex in
$\P$ is the whole $\P$). However, in the process of barycentric
subdivision of $\P$ we first make the stellar subdivisions at
2-dimensional simplices, and the stars of 1-simplices in the
resulting complex are simplicial complexes. Note also that if
$\st_{\P}\sigma$ is not a simplicial complex, then the map
$\beta\colon\k[\P]\to\k[\P']$ is not determined by the conditions
specified in Lemma~\ref{lemm:frppt}. That is, the images of
degree-two generators do not determine the images of generators of
higher degree. However, it is still possible to define the map
$\beta\colon\k[\P]\to\k[\P']$ for an arbitrary simplicial poset;
see Section~\ref{blowu}.
\end{exam}

\begin{lemm}\label{lemm:cmtil}
Assume that $\k[\P]$ is a Cohen--Macaulay ring, and $\widetilde\P$
a stellar subdivision of $\P$ at $\sigma$ such that $\st_\P\sigma$
is a simplicial complex. Then
\begin{itemize}
\item[(a)] $\st_\P\sigma$ is a Cohen--Macaulay complex;

\item[(b)] $\k[\til\P]$ is a Cohen--Macaulay ring.
\end{itemize}
\end{lemm}
\begin{proof}
(a) As $\st_\P\sigma=\sigma\mathbin{*}\lk_\P\sigma$, it is enough
to prove that $\lk_\P\sigma$ is Cohen--Macaulay. This can be done
by showing that the simplicial homology of $\lk_\P\sigma$ is a
direct summand in the local cohomology of $\k[\P]$, as in the
proof of Hochster's theorem (see~\cite[Th.~II.4.1]{stan96}
or~\cite[Th.~5.3.8]{br-he98}).

(b) (Compare proof of Lemma 9.2 of~\cite{ma-pa06}.) Choose an lsop
$\theta_1,\ldots,\theta_n\in \k[\P]$ (we can always assume that an
lsop exists by passing to an infinite extension field,
see~\cite[Th.~2.1.10]{br-he98}) and denote
$\tilde\theta_i=\beta(\theta_i)$. Applying the functors
$\otimes_{\k[\theta_1,\ldots,\theta_n]}\k$ and
$\otimes_{\k[\tilde\theta_1,\ldots,\tilde\theta_n]}\k$ to
diagram~\eqref{exrow}, we get a map between the long exact
sequences for $\Tor$. Consider the following fragment (we denote
$\Tor_\theta=\Tor_{\k[\theta_1,\ldots,\theta_n]}$):
\[
\begin{array}{cccc}
  \Tor^{-2}_\theta(\k[\st\sigma],\k)
  \stackrel{f}{\to} &\!\Tor^{-1}_\theta(R,\k)
  \to &\! \Tor^{-1}_\theta(\k[\P],\k)  \to &\!
  \Tor^{-1}_\theta(\k[\st\sigma],\k)\\
  \downarrow & \downarrow{\scriptstyle\cong} & \downarrow & \downarrow\\
  \Tor^{-2}_{\tilde\theta}(\k[\st v],\k)
  \stackrel{\tilde f}{\to} &\!
  \Tor^{-1}_{\tilde\theta}(\widetilde R,\k)  \to &\!
  \Tor^{-1}_{\tilde\theta}(\k[\widetilde\P],\k)  \to &\!
  \Tor^{-1}_{\tilde\theta}(\k[\st v],\k).
\end{array}
\]
Since $\k[\P]$ is Cohen--Macaulay,
$\Tor^{-1}_\theta(\k[\P],\k)=0$, therefore, $f$ is onto. Then
$\tilde f$ is also onto. Since $\st_{\P}\sigma$ (and
$\st_{\widetilde\P}v$) is a simplicial complex and
$|\st_{\P}\sigma|\cong|\st_{\widetilde\P}v|$, part~(a) of this
lemma and Theorem~\ref{theo:munkr} imply that $\k[\st v]$ is
Cohen--Macaulay. Therefore, $\Tor^{-1}_{\tilde\theta}(\k[\st
v],\k)=0$. Since $\tilde f$ is surjective, we also have
$\Tor^{-1}_{\tilde\theta}(\k[\widetilde\P],\k)=0$. Then
$\k[\widetilde\P]$ is free as a module over
$\k[\tilde\theta_1,\ldots,\tilde\theta_n]$
(see~\cite[Lemma~VII.6.2]{macl63}) and thereby is Cohen--Macaulay.
\end{proof}

\begin{theo}\label{theo:cmpos}
The face ring $\k[\P]$ of a simplicial poset $\P$ is
Cohen--\-Ma\-ca\-u\-lay if and only if $\P$ is Cohen--Macaulay.
\end{theo}
\begin{proof}
Assume $\k[\P]$ is Cohen--Macaulay. Since the barycentric
subdivision $\P'$ is obtained by a sequence of stellar
subdivisions, subsequent application of Lemma~\ref{lemm:cmtil}
gives that $\k[\P']$ is also Cohen--Macaulay. Hence, $\P$ is a
Cohen--Macaulay poset. The converse statement
is~\cite[Cor.~3.7]{stan91}.
\end{proof}

\section{Pseudomanifolds and orientations}\label{spseu}
The question of identifying the class of simplicial posets which
arise as $\mathcal P(\Gamma)$ for some torus graph $\Gamma$ might
be a difficult one, although our next statement sheds some light
on this problem. The following is a straightforward extension of
the notion of \emph{pseudomanifold}~\cite[Def.~0.3.15]{stan96} to
simplicial posets.

\begin{defi}
A simplicial poset $\P$ is called an \emph{$(n-1)$-dimensional
pseudomanifold} (without boundary) if
\begin{itemize}
\item[(a)] for every element $\sigma\in\P$ there is an element
$\tau$ of rank $n$ such that $\sigma\le\tau$ (in other words, $\P$
is \emph{pure $(n-1)$-dimensional});

\item[(b)] for every element $\sigma\in\P$ of rank $(n-1)$ there
are exactly two elements $\tau$ of rank $n$ such that
$\sigma<\tau$;

\item[(c)] for every two elements $\tau$ and $\tau'$ of rank $n$
there is a sequence $\tau=\tau_1,\tau_2,\ldots,\tau_k=\tau'$ of
elements such that $\rank\tau_i=n$ and $\tau_i\wedge\tau_{i+1}$
contains an element of rank $(n-1)$ for $i=1,\ldots,k-1$.
\end{itemize}
\end{defi}

Examples of pseudomanifolds are provided by triangulations or
simplicial cell decompositions of topological manifolds. However,
not every pseudomanifold arises in this way, see
Example~\ref{expsm} below.

\begin{theo}\label{theo:eqcal}
\begin{itemize}
\item[(a)] Let $\Gamma$ be a torus graph; then $\P(\Gamma)$ is a
pseu\-do\-ma\-ni\-fold, and the face ring $\Z(\P)$ admits an lsop;

\item[(b)] Given an arbitrary pseudomanifold $\P$ and an lsop in
$\Z(\P)$, one can canonically construct a torus graph $\Gamma_\P$.
\end{itemize}
Moreover, $\Gamma_{\P(\Gamma)}=\Gamma$.
\end{theo}
\begin{proof}
(a) Vertices of $\P(\Gamma)$ correspond to $(n-1)$-faces
of~$\Gamma$. As every face of $\Gamma$ contains a vertex and
$\Gamma$ is $n$-valent, $\P(\Gamma)$ is pure $(n-1)$-dimensional.
Condition~(b) from the definition of a pseudomanifold follows from
the fact that every edge of $\Gamma$ has exactly two vertices,
while~(c) follows from the connectivity of~$\Gamma$. In order to
find an lsop, we identify $\Z[\P(\Gamma)]$ with a subset of
$H^*(BT)^{V(\Gamma)}$ (see Theorem~\ref{tequi}) and consider the
constant map $c\colon H^*(BT)\to H^*(BT)^{V(\Gamma)}$. It factors
through a monomorphism $H^*(BT)\to \Z[\P(\Gamma)]$, and
Proposition~\ref{prop:dlsop} guarantees that the $c$-image of a
basis in $H^*(BT)$ is an lsop.

(b) Let $\P$ be a pseudomanifold of dimension $(n-1)$. Define a
graph $\Gamma_\P$ whose vertices correspond to $(n-1)$-dimensional
simplices $\sigma\in\P$, and in which the number of edges between
two vertices $\sigma$ and $\sigma'$ equals the number of
$(n-2)$-dimensional simplices in $\sigma\wedge\sigma'$. Then $\Gamma_\P$ is
a connected $n$-valent graph, and we need to define an axial
function. The following argument is similar to that
of~\cite[\S3]{masu05}, compare also a similar treatment of ``edge
vectors'' in~\cite{pano01}. We can regard an lsop as a map
$\lambda\colon H^*(BT)\to\Z[\P]$. As usual, assume that $\P$ has
$m$ vertices and let $v_1,\ldots,v_m$ be the corresponding
degree-two generators of $\Z[\P]$. Then for $t\in H^2(BT)$ we can
write
\[
  \lambda(t)=\sum_{i=1}^m\lambda_i(t)v_i,
\]
where $\lambda_i$ is a linear function on $H^2(BT)$, that is, an
element of $H_2(BT)$. Let $e$ be an oriented edge of $\Gamma$ and
$p=i(e)$ its initial vertex. This vertex corresponds to an
$(n-1)$-simplex of $\P$, and we denote by
$I(p)\subset\{1,\ldots,m\}$ the set of its vertices in~$\P$; note
that $|I(p)|=n$. Since $\lambda$ is an lsop, the set
$\{\lambda_i\colon i\in I(p)\}$ is a basis in $H_2(BT)$. Now we
define the axial function $\alpha\colon E(\Gamma)\to H^2(BT)$ by
requiring that its value on $E(\Gamma)_p$ is the dual basis to
$\{\lambda_i\colon i\in I(p)\}$. In more detail, the edge $e$
corresponds to an $(n-2)$-simplex of $\P$ and let $\ell\in I(p)$
be the unique vertex which is not in this $(n-2)$-simplex. Then we
define $\alpha(e)$ by requiring that
\begin{equation}\label{eaxfu}
  \langle\alpha(e),\lambda_i\rangle=\delta_{i\ell},\quad i\in
  I(p),
\end{equation}
where $\delta_{i\ell}$ is the Kronecker delta. Now we have to
check the three conditions from the definition of a torus axial
function. Let $p'=t(e)=i(\bar e)$. Note that the intersection of
$I(p)$ and $I(p')$ consists of at least $(n-1)$ elements. If
$I(p)=I(p')$ then $\Gamma$ has only two vertices, like in
Example~\ref{ebana}, while $\P$ is obtained by gluing together two
$(n-1)$-simplices along their boundaries, see Example~\ref{2tria}.
Otherwise, $|I(p)\cap I(p')|=n-1$ and we have $\ell\notin I(p')$.
Let $\ell'$ be an element such that $\ell'\in I(p')$ but
$\ell'\notin I(p)$. Then~\eqref{eaxfu} guarantees that
$\langle\alpha(e),\lambda_i\rangle= \langle\alpha(\bar
e),\lambda_i\rangle=0$ for $i\in I(p)\cap I(p')$. As we work with
integral bases, this implies $\alpha(\bar e)=\pm\alpha(e)$. It
also follows that $\alpha(E(\Gamma)_p\setminus e)$ and
$\alpha(E(\Gamma)_{p'}\setminus\bar e)$ give the same bases in the
quotient space $H^2(BT)/\alpha(e)$. Identifying these bases, we
obtain a connection $\theta_e\colon E(\Gamma)_p\to E(\Gamma)_{p'}$
satisfying $\alpha(\theta_e(e'))\equiv\alpha(e')\mod\alpha(e)$ for
any $e'\in E(\Gamma)_p$, as needed. The rest of the statement is
straightforward.
\end{proof}

Note that the above theorem does not give a complete
characterisation of simplicial posets of the form $\mathcal
P(\Gamma)$, as it may happen that $\P(\Gamma_\P)\ne\P$. In fact,
here is a counterexample.

\begin{exam}\label{expsm}
Let $\P$ be a triangulation of a 2-dimensional sphere different
from the boundary of a simplex. Choose two vertices that are not
joined by an edge. Let $\P'$ be the complex obtained by
identifying these two vertices. Then $\P'$ is a pseudomanifold. If
$\Z[\P]$ admits an lsop, so does $\Z[\P']$ (this easily follows
from Proposition~\ref{prop:dlsop}). However,
$\P(\Gamma_{\P'})\ne\P'$ (in fact, $\P(\Gamma_{\P'})=\P$). It
follows that $\P'$ does not arise from any torus graph.
\end{exam}

\begin{defi}
We say that an assignment $\o\colon V(\G)\to \{\pm 1\}$ is an
\emph{orientation} on $\mG$ if $\o(i(e))\alpha(e)=-\o(i(\bar
e))\alpha(\bar e)$ for every $e\in E(\G)$.
\end{defi}

\begin{exam}
Let $M$ be a torus manifold which admits a $T$-invariant almost
complex structure. The almost complex structure induce
orientations on $M$ and its characteristic submanifolds (an
omniorientation). The associated torus axial function $\alpha_M$
satisfies $\alpha_M(\bar e)=-\alpha_M(e)$ for any oriented edge
$e$. In this case we can take $o(p)=1$ for every $p\in V(\G_M)$.
\end{exam}

\begin{prop}
An omniorientation of a torus manifold $M$ induces an orientation
of the associated torus graph~$\Gamma_M$.
\end{prop}
\begin{proof}
Given a vertex $p\in M^T=V(\Gamma_M)$ we set $o(p)=1$ if the
canonical orientation of the sum of complex one-dimensional
representation spaces in the right hand side of~\eqref{eqn:deco}
coincides with the orientation of $\tau_p M$ induced by the
orientation of $M$, and set $o(p)=-1$ otherwise.
\end{proof}

\begin{exam}
Let $\Gamma$ be a complete graph on four vertices
$p_1,p_2,p_3,p_4$. Choose a basis $t_1,t_2,t_3\in H^2(BT^3)$ and
define an axial function by setting
\begin{gather*}
  \alpha(p_1p_2)=\alpha(p_3p_4)=t_1, \
  \alpha(p_1p_3)=\alpha(p_2p_4)=t_2, \
  \alpha(p_1p_4)=\alpha(p_2p_3)=t_3
\end{gather*}
and $\alpha(\bar e)=\alpha(e)$ for any oriented edge~$e$. A direct
check shows that this torus graph is not orientable. In fact, this
graph is associated with the pseudomanifold (simplicial cell
complex) shown on Fig.~\ref{frp2c} via the construction of
Theorem~\ref{theo:eqcal}~(b). This pseudomanifold $\mathcal P$ is
homeomorphic to $\R P^2$ (the opposite outer edges are identified
according to the arrows shown), the ring $\Z[\P]$ has three
two-dimensional generators $v_p,v_q,v_r$, which constitute an
lsop. Note that $\R P^2$ itself is non-orientable; in fact, this
example is generalised by the following Proposition.
\begin{figure}[h]
\begin{picture}(120,30)
\put(45,0){\vector(1,0){29.5}}
\put(75,0){\vector(0,1){29.5}}
\put(75,30){\vector(-1,0){29.5}}
\put(45,30){\vector(0,-1){29.5}}
\put(45,0){\line(1,1){30}}
\put(45,30){\line(1,-1){30}}
\put(45,0){\circle*{1}}
\put(45,30){\circle*{1}}
\put(75,0){\circle*{1}}
\put(75,30){\circle*{1}}
\put(60,15){\circle*{1}}
\put(42,-1){$q$}
\put(76,-1){$p$}
\put(42,29){$p$}
\put(76,29){$q$}
\put(60,16){$r$}
\end{picture}
\caption{Simplicial cell decomposition of $\mathbb RP^2$ with 3
vertices.}
\label{frp2c}
\end{figure}
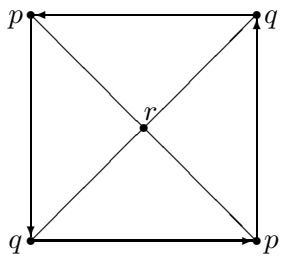
\end{exam}

\begin{prop} A torus graph
$\Gamma$ is orientable if and only if the associated
pseudomanifold $\mathcal P(\Gamma)$ is orientable.
\end{prop}
\begin{proof}
Let $p\in V(\Gamma)$ and $\sigma$ the corresponding
$(n-1)$-simplex of $\mathcal P(\Gamma)$. There is a canonical
one-to-one correspondence between $E(\Gamma)_p$ and the vertex set
of $\sigma$, see the proof of Theorem~\ref{theo:eqcal}~(b). Choose
a basis of $H^2(BT)$. Assume first that $\mathcal P(\Gamma)$ is
oriented. Choose a ``positive'' (that is, compatible with the
orientation) order of vertices of $\sigma$; this allows to regard
$\alpha(E(\Gamma)_p)$ as a basis of $H^2(BT)$. We set $o(p)=1$ if
this is a positively oriented basis, and $o(p)=-1$ otherwise. This
defines an orientation on $\Gamma$. To prove the opposite
direction we just reverse this procedure.
\end{proof}

\section{Blow-ups of torus manifolds and torus graphs}\label{blowu}
Here we relate the following three geometric constructions:
\begin{itemize}
\item[(a)] blowing up a torus manifold at a facial
submanifold~\cite[\S9]{ma-pa06};

\item[(b)] cutting a face from a simple polytope or, more
generally, blowing up a GKM graph~\cite[\S2.2]{gu-za99};

\item[(c)] stellar subdivision of a simplicial poset.
\end{itemize}

Let $M$ be a torus manifold with the projection map $\pi\colon
M\to Q$ onto the orbit space, and $F$ a face of $Q$ (details may
be found in~\cite{ha-ma03} or~\cite{ma-pa06}; a reader less
familiar with torus manifolds may assume $M$ to be a smooth
projective toric variety, in which case $Q$ is a convex simple
polytope). Replacing the facial submanifold $M_F=\pi^{-1}(F)$ of
$M$ by the complex projectivisation $P(\n{F})$ of its normal
bundle $\n{F}$, we obtain a new torus manifold $\widetilde M$. The
passage from $M$ to $\widetilde M$ is called \emph{blowing-up} of
$M$ at $M_F$. The orbit space $\widetilde Q$ of $\widetilde M$ is
then obtained by ``cutting off" the face $F$ from $Q$. As
explained in Example~\ref{etmtg}, the 1-skeleton of $Q$ is a torus
graph. The general construction of blow-up of a GKM-graph is
described in~\cite[\S2.2.1]{gu-za99}; in particular, it applies to
torus graphs and agrees with the topological picture for the
graphs coming from manifolds. We briefly review their construction
below, and illustrate it in a couple of examples.

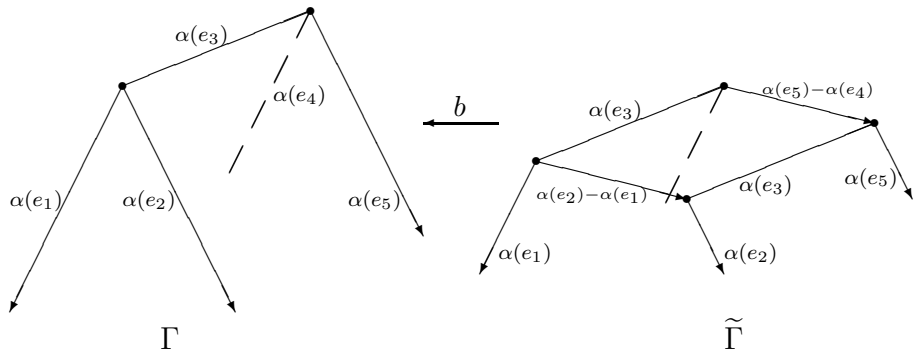
\begin{figure}[h]
\begin{picture}(120,40)
\put(15,30){\vector(-1,-2){15}}
\put(15,30){\vector(1,-2){15}}
\put(15,30){\line(5,2){25}}
\put(40,40){\vector(1,-2){15}}
\multiput(40,40)(-3,-6){4}{\line(-1,-2){1.8}}
\put(15,30){\circle*{1}}
\put(40,40){\circle*{1}}
\put(70,20){\vector(-1,-2){7.5}}
\put(70,20){\vector(4,-1){20}}
\put(90,15){\vector(1,-2){5}}
\put(70,20){\line(5,2){25}}
\put(90,15){\line(5,2){25}}
\put(95,30){\vector(4,-1){20}}
\put(115,25){\vector(1,-2){5}}
\multiput(95,30)(-3,-6){3}{\line(-1,-2){1.8}}
\put(70,20){\circle*{1}}
\put(90,15){\circle*{1}}
\put(95,30){\circle*{1}}
\put(115,25){\circle*{1}}
\put(0,14){$\scriptstyle\alpha(e_1)$}
\put(15,14){$\scriptstyle\alpha(e_2)$}
\put(22,36){$\scriptstyle\alpha(e_3)$}
\put(35,28){$\scriptstyle\alpha(e_4)$}
\put(45,14){$\scriptstyle\alpha(e_5)$}
\put(65,7){$\scriptstyle\alpha(e_1)$}
\put(95,7){$\scriptstyle\alpha(e_2)$}
\put(70,15){$\scriptscriptstyle\alpha(e_2)-\alpha(e_1)$}
\put(97,16){$\scriptstyle\alpha(e_3)$}
\put(77,26){$\scriptstyle\alpha(e_3)$}
\put(111,17){$\scriptstyle\alpha(e_5)$}
\put(100,29){$\scriptscriptstyle\alpha(e_5)-\alpha(e_4)$}
\put(65,25){\vector(-1,0){10}} \put(59,26){$b$}
\put(20,-5){\large$\Gamma$}
\put(95,-5){\large$\widetilde\Gamma$}
\end{picture}
\medskip
\caption{Blow up at an edge} \label{fblo1}
\end{figure}

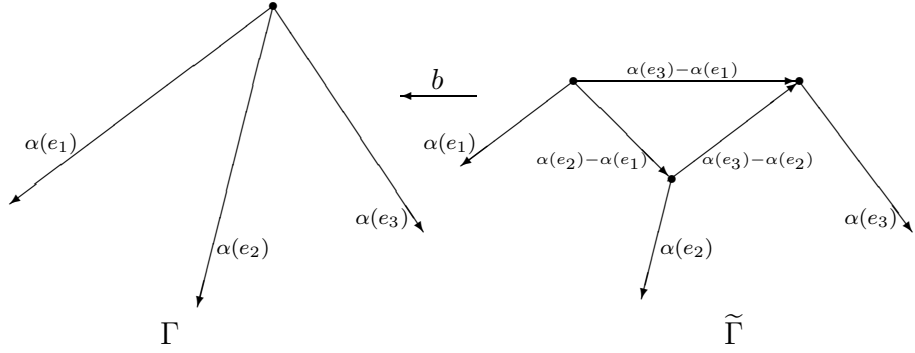
\begin{figure}[h]
\begin{picture}(120,40)
\put(35,40){\vector(-4,-3){35}}
\put(35,40){\vector(-1,-4){10}}
\put(35,40){\vector(2,-3){20}}
\put(35,40){\circle*{1}}
\put(75,30){\vector(-4,-3){15}}
\put(75,30){\vector(1,-1){12.6}}
\put(75,30){\vector(1,0){29.5}}
\put(88,17){\vector(-1,-4){4}}
\put(87.5,16.8){\vector(4,3){17.2}}
\put(105,30){\vector(3,-4){15}}
\put(75,30){\circle*{1}}
\put(88,17){\circle*{1}}
\put(105,30){\circle*{1}}
\put(2,21){$\scriptstyle\alpha(e_1)$}
\put(27.5,7){$\scriptstyle\alpha(e_2)$}
\put(46,11){$\scriptstyle\alpha(e_3)$}
\put(55,21){$\scriptstyle\alpha(e_1)$}
\put(86.5,7){$\scriptstyle\alpha(e_2)$}
\put(70,19){$\scriptscriptstyle\alpha(e_2)-\alpha(e_1)$}
\put(111,11){$\scriptstyle\alpha(e_3)$}
\put(82,31){$\scriptscriptstyle\alpha(e_3)-\alpha(e_1)$}
\put(92,19){$\scriptscriptstyle\alpha(e_3)-\alpha(e_2)$}
\put(62,28){\vector(-1,0){10}}
\put(56,29){$b$}
\put(20,-5){\large$\Gamma$}
\put(95,-5){\large$\widetilde\Gamma$}
\end{picture}
\medskip
\caption{Blow up at a vertex} \label{fblo2}
\end{figure}

Let $F$ be a $k$-face of $\Gamma$ (of codimension $n-k$). The
\emph{blow-up} of $\Gamma$ at $F$, denoted $\widetilde\Gamma$, has
vertex set $V(\widetilde\Gamma)=(V(\Gamma)\setminus V(F))\cup
V(F)^{n-k}$, that is, each vertex $p\in V(F)$ is replaced by
$(n-k)$ vertices $\til p_1,\ldots,\til p_{n-k}$. It is convenient
to regard those points as chosen close to $p$ on edges from
$E_p(\Gamma)\setminus E_p(F)$, and we denote by $p'_i$ the
endpoint of the edge containing both $p$ and~$\til p_i$, \
$i=1,\ldots,n-k$. (We also assume $\theta_{pq}(pp'_i)=qq'_i$ if
$p$ and $q$ are joined by an edge in~$F$.) Then we have four types
of edges in $\til\Gamma$, and the corresponding values of the
axial function $\til\alpha\colon E(\til\Gamma)\to H^*(BT)$:
\begin{itemize}
\item[(a)] $\til p_i\til p_j$ for every $p\in V(F)$; \
$\til\alpha(\til p_i\til p_j)=\alpha(pp'_j)-\alpha(pp'_i)$;

\item[(b)] $\til p_i\til q_i$ if $p$ and $q$ were joined by an
edge in~$F$; \ $\til\alpha(\til p_i\til q_i)=\alpha(pq)$;

\item[(c)] $\til p_ip'_i$ for every $p\in V(F)$; \
$\til\alpha(\til p_ip'_i)=\alpha(pp'_i)$;

\item[(d)] edges ``coming from $\Gamma$'', that is, $e\in
E(\Gamma)$ such that $i(e)\notin V(F)$ and $t(e)\notin V(F)$; \
$\til\alpha(e)=\alpha(e)$,
\end{itemize}
see Fig.~\ref{fblo1} ($n=3$, $k=1$) and Fig.~\ref{fblo2} ($n=3$,
$k=0$).

There is a \emph{blow-down map} $b\colon\til\Gamma\to\Gamma$
preserving the face structure. The face $F\subset\Gamma$ is blown
up to a new facet $\til F\subset\til\Gamma$ (unless $F$ itself was
a facet, in which case $\til\Gamma=\Gamma$). For every face
$H\subset\Gamma$ which is not contained in $F$, there is a unique
face $\til H\subset\til\Gamma$ that is mapped onto~$H$. The
blow-down map induces an equivariant cohomology map $b^*\colon
H^*_T(\Gamma)\to H^*_T(\til\Gamma)$. In fact, this map can be
easily identified by the following commutative diagram
\begin{equation}
\label{egrre}
\begin{CD}
  H^*_T(\Gamma) @>b^*>> H^*_T(\widetilde\Gamma)\\
  @VrVV @VV{\til r}V\\
  H^*(BT)^{V(\Gamma)} @>V(b)^*>> H^*(BT)^{V(\til\Gamma)}
\end{CD}
\end{equation}
(compare~\eqref{epore}), where $r$ and $\til r$ are the
monomorphisms from the definition of equivariant cohomology of a
torus graph, and $V(b)^*$ is the map induced by the set map
$V(b)\colon V(\til\Gamma)\to V(\Gamma)$. The next lemma describes
the images of the two-dimensional generators $\ta{G}\in
H^*_T(\Gamma)$ corresponding to the facets $G\subset\Gamma$.

\begin{lemm}
Given a facet $G\subset\Gamma$, we have $b^*(\ta{G})=\ta{\til
F}+\ta{\til G}$ if $F\subset G$ and $b^*(\ta{G})=\ta{\til G}$
otherwise.
\end{lemm}
\begin{proof}
We use~\eqref{egrre} and check that the images of $\ta{G}$ and
$\ta{\til F}+\ta{\til G}$ (or $\ta{\til G}$) under the horizontal
maps agree. Let $p\in V(\Gamma)$ be a vertex. If $p\notin F$, then
$b^{-1}(p)=p$ and $r(\ta{G})(p)=\til r(\ta{\til G})(p)$, \ $\til
r(\ta{\til F})(p)=0$. Thus we may assume $p\in F$, and then we have
$b(\til p_i)=p$, \ $i=1,\ldots,n-k$.

First consider the case $F\not\subset G$. If $p\notin G$, then
$r(\ta{G})(p)=\til r(\ta{\til G})(\til p_i)=0$. Otherwise $p\in
G\cap F$. Let $e$ be the unique edge such that $e\in E_p(\Gamma)$
and $e\notin G$. Then $e=pq$ for some $q\in V(F)$ (because
$F\not\subset G$). From~\eqref{etaf1} we obtain
\[
  r(\ta{G})(p)=\alpha(pq)=\til\alpha(\til p_i\til q_i)=
  \til r(\ta{\til G})(\til p_i),
\]
see Fig.~\ref{fres1}. It follows that $V(b)^* r(\ta{G})=\til
r(\ta{\til G})$, and therefore, $b^*(\ta{G})=\ta{\til G}$.
\begin{figure}[h]
\begin{picture}(120,40)
\put(15,30){\line(-1,-2){15}}
\put(15,30){\line(1,-2){15}}
\put(15,30){\line(5,2){25}}
\put(40,40){\line(1,-2){15}}
\put(15,30){\circle*{1}}
\put(40,40){\circle*{1}}
\put(70,20){\line(-1,-2){7.5}}
\put(70,20){\line(4,-1){20}}
\put(90,15){\line(1,-2){5}}
\put(70,20){\line(5,2){25}}
\put(90,15){\line(5,2){25}}
\put(95,30){\line(4,-1){20}}
\put(115,25){\line(1,-2){5}}
\put(70,20){\circle*{1}}
\put(90,15){\circle*{1}}
\put(95,30){\circle*{1}}
\put(115,25){\circle*{1}}
\put(26,35.5){$F$}
\put(13,31){$p$}
\put(38.5,37){$q$}
\put(14,10){$G$}
\put(92,21){$\til F$}
\put(88,11){$\til p_i$}
\put(113,21){$\til q_i$}
\put(77,10){$\til G$}
\put(65,25){\vector(-1,0){10}}
\put(59,26){$b$}
\end{picture}
\caption{}
\label{fres1}
\end{figure}

Now let $F\subset G$. In this case the unique edge $e$ such that
$e\in E_p(\Gamma)$ and $e\notin G$ is of type $pp'_j$, see
Fig.~\ref{fres2}. Using~\eqref{etaf1} we calculate
\begin{align*}
  r(\ta{G})(p)&=\alpha(pp'_j),\\
  \til r(\ta{\til G})(\til p_i)&=\til\alpha(\til p_i\til p_j)
  =\alpha(pp'_j)-\alpha(pp'_i),\\
  \til r(\ta{\til F})(\til p_i)&=\til\alpha(\til
  p_ip'_i)=\alpha(pp'_i).
\end{align*}
As in the previous case it follows that $V(b)^* r(\ta{G})=\til
r(\ta{\til F})+\til r(\ta{\til G})$, and therefore,
$b^*(\ta{G})=\ta{\til F}+\ta{\til G}$.
\begin{figure}[h]
\begin{picture}(120,40)
\put(15,30){\line(-1,-2){15}}
\put(15,30){\line(1,-2){15}}
\put(15,30){\line(5,2){25}}
\put(40,40){\line(1,-2){15}}
\put(15,30){\circle*{1}}
\put(40,40){\circle*{1}}
\put(0,0){\circle*{1}}
\put(30,0){\circle*{1}}
\put(70,20){\vector(-1,-2){7.3}}
\put(70,20){\line(4,-1){20}}
\put(90,15){\vector(1,-2){4.8}}
\put(70,20){\line(5,2){25}}
\put(90,15){\line(5,2){25}}
\put(90,15){\vector(-4,1){19.5}}
\put(115,25){\vector(-4,1){19.5}}
\put(115,25){\line(1,-2){5}}
\put(70,20){\circle*{1}}
\put(90,15){\circle*{1}}
\put(95,30){\circle*{1}}
\put(115,25){\circle*{1}}
\put(62.5,5){\circle*{1}}
\put(95,5){\circle*{1}}
\put(26,35.5){$F$}
\put(13,31){$p$}
\put(1.5,0){$p'_j$}
\put(25.5,0){$p'_i$}
\put(35,15){$G$}
\put(92,21){$\til F$}
\put(88,11){$\til p_i$}
\put(70,16){$\til p_j$}
\put(91,4){$p'_i$}
\put(64,4){$p'_j$}
\put(105,12){$\til G$}
\put(65,25){\vector(-1,0){10}}
\put(59,26){$b$}
\end{picture}
\caption{} \label{fres2}
\end{figure}
\end{proof}

\begin{coro}
After the identifications $H^*_T(\Gamma)\cong\Z[\mathcal P(\Gamma)]$ and
$H^*_T(\til\Gamma)\cong\Z[\mathcal P(\til\Gamma)]$, the equivariant
cohomology map induced by the blow-down
$b\colon\widetilde\Gamma\to\Gamma$ coincides with the map $\beta$
from Lemma~\ref{lemm:frppt}.
\end{coro}
\begin{proof}
Remember that the poset $\mathcal P(\Gamma)$ is formed by the
faces of $\Gamma$ with the reversed inclusion relation, and the
isomorphism $H^*_T(\Gamma)\cong\Z[\mathcal P(\Gamma)]$ is
established by identifying $\ta{H}$ with $\v{H}$ for all faces
$H\subset\Gamma$. Let $\sigma\in\mathcal P(\Gamma)$ be the element
corresponding to the face~$F$. Then an element $\tau\in\mathcal
P(\Gamma)$ satisfies $\tau\in\st_{\P}\sigma$ if and only if the
corresponding face $H\subset\Gamma$ satisfies $F\cap
H\ne\varnothing$. The degree-two generators $v_i$, \
$i=1,\ldots,m$, of $\Z[\P(\Gamma)]$ (or $\Z[\P(\til\Gamma)]$)
correspond to the generators $\ta{G}$ of $H^*_T(\Gamma)$ (or
$\ta{\til G}$ of $H^*_T(\til\Gamma)$ respectively). Making the
appropriate identifications, we see that the map from
Lemma~\ref{lemm:frppt} is determined by the conditions
\begin{align*}
  \ta{H} & \mapsto\ta{H} &&\text{if \ }F\cap H=\varnothing,\\
  \ta{G} & \mapsto\ta{\til F}+\ta{\til G} &&\text{if \ }F\subset G,\\
  \ta{G} & \mapsto\ta{\til G} &&\text{if \ }F\not\subset G.
\end{align*}
The blow-down map $b^*$ satisfies these conditions, whence the
proof follows.
\end{proof}

{\sloppy Returning to torus manifolds,
in~\cite[Lemma~9.2]{ma-pa06} we proved that $H^{odd}(M)=0$ implies
$H^{odd}(\widetilde M)=0$. Also, $H^{odd}(M)=0$ implies that
$\Z[\P]$ is Cohen--Macaulay~\cite[Th.~7.7]{ma-pa06} (here $\P$ is
the face poset of the orbit space~$Q$). Now the above mentioned
analogy between the proof of~\cite[Lemma~9.2]{ma-pa06} and the
proof of Lemma~\ref{lemm:cmtil} becomes even more transparent.
Note that we have isomorphisms $H^*_T(M_F)\cong\Z[\st_{\P}\sigma]$
and, $H^*_T(\widetilde M_{\widetilde F})\cong\Z[\st_{\til\P}v]$.
We are also ready to give the proof of the other direction
of~\cite[Lemma~9.2]{ma-pa06}, promised in the end of Section~9
of~\cite{ma-pa06}.

}
\begin{lemm}
$H^{odd}(\widetilde M)=0$ if and only if $H^{odd}(M)=0$.
\end{lemm}
\begin{proof}
Assume $H^{odd}(\widetilde M)=0$. Then $\Z[\widetilde\P]$ is
Cohen--Macaulay by Theorem~7.7 of~\cite{ma-pa06}. We claim that
$\Z[\P]$ is also Cohen-Macaulay (i.e, the converse of
Lemma~\ref{lemm:cmtil} holds). Indeed, by
Theorem~\ref{theo:cmpos}, $\widetilde\P$ is a Cohen--Macaulay
poset. Choose a simplicial complex $\mathcal S$ which is a common
subdivision of $\widetilde\P$ and $\P$ (for example, we can take
$\mathcal S$ to be the barycentric subdivision of $\widetilde\P$).
By Corollary~\ref{coro:cmcha}, $\mathcal S$ is a Cohen--Macaulay
complex, whence $\P$ is a Cohen--Macaulay poset. Applying
Theorem~\ref{theo:cmpos} again we get that $\Z[\P]$ is
Cohen--Macaulay. Then $H^{odd}(M)=0$ by Theorem~7.7
of~\cite{ma-pa06}. The other direction of the lemma is already
proven in~\cite{ma-pa06}.
\end{proof}

\section{Dehn-Sommerville equations}
Let $\P$ be a simplicial poset of rank $n$ (i.e. of dimension
$n-1$). Let $f_i$ denote the number of $i$-dimensional simplices
in $\P$, \ $0\le i\le n-1$. Since $\P$ has a unique initial
element $\hatzero$, we have $f_{-1}=1$. The \emph{$h$-vector} $\mb
h(\P)=(h_0,\ldots,h_n)$ of $\P$ is defined from the polynomial
identity
\begin{equation} \label{eqn:h-vector}
  \sum_{i=0}^nh_it^{n-i}=\sum_{i=0}^n f_{i-1}(t-1)^{n-i}.
\end{equation}

Let $\P_{\ge\sigma}=\{\tau\in\P\colon\tau\ge\sigma\}$ be the
subposet $\P$ with the induced rank function. For a simplex
$\sigma\in\P$ we set
\begin{equation} \label{eqn:tchi}
\chi(\P_{\ge\sigma}):=\sum_{\tau\ge \sigma}(-1)^{\rank \tau-1}.
\end{equation}

\begin{theo} \label{theo:DS}
$\displaystyle{\sum_{i=0}^n(h_{n-i}-h_{i})t^i=\sum_{\sigma\in\P}
\Bigl(1+(-1)^n\chi(\P_{\ge\sigma})\Bigr)(t-1)^{n-\rank \sigma}}$.
\newline
In particular, the Dehn-Sommerville equations $h_i=h_{n-i}$ hold
if $\chi(\P_{\ge\sigma})=(-1)^{n-1}$ for every $\sigma\in\P$.
\end{theo}
\begin{proof}The argument below is essentially the same as that used by
Hibi in~\cite[p.~91]{hibi95}. We have
\begin{equation} \label{eqn:DS}
\begin{split}
\sum_{i=0}^n h_it^i&=t^n \sum_{i=0}^n h_i(1/t)^{n-i}
=t^n\sum_{i=0}^n f_{i-1}((1-t)/t)^{n-i}\qquad \text{by (\ref{eqn:h-vector})}\\
&=\sum_{i=0}^n f_{i-1}t^i(1-t)^{n-i}
=\sum_{\tau\in\P} t^{\rank \tau}(1-t)^{n-\rank \tau}\\
&=\sum_{\tau\in\P}\sum_{\sigma\le \tau}(t-1)^{\rank \tau-\rank
\sigma}
(1-t)^{n-\rank \tau}\\
&=\sum_{\tau\in\P}\sum_{\sigma\le \tau}(-1)^{n-\rank
\tau}(t-1)^{n-\rank \sigma}
\\
&=\sum_{\sigma\in\P}(t-1)^{n-\rank \sigma}\sum_{\tau\ge \sigma}
(-1)^{n-\rank \tau}\\
&=\sum_{\sigma\in\P}(t-1)^{n-\rank
\sigma}(-1)^{n-1}\chi(\P_{\ge\sigma}) \qquad \text{by
(\ref{eqn:tchi})}
\end{split}
\end{equation}
where the fifth equality follows from the binomial expansion of
the right hand side of the identity
$t^{\rank\tau}=((t-1)+1)^{\rank\tau}$.

On the other hand, we have
\begin{equation} \label{eqn:DS2}
\sum_{i=0}^nh_{n-i}t^i=\sum_{i=0}^n
h_it^{n-i}=\sum_{i=0}^nf_{i-1}(t-1)^{n-i}
=\sum_{\sigma\in\P}(t-1)^{n-\rank \sigma}.
\end{equation}

Subtracting (\ref{eqn:DS}) from (\ref{eqn:DS2}) we obtain the
theorem.
\end{proof}

\begin{coro}[\cite{bu-pa02}]
If $K$ is a triangulation of a closed $(n-1)$-manifold, then
\[
h_{n-i}-h_i=(-1)^i\binom{n}{i}\big(\chi(K)-\chi(S^{n-1})\big).
\]
\end{coro}
\begin{proof}
Let $\P$ be the face poset of $K$ with an added initial element
(corresponding to the empty simplex). Then for any $\sigma\in\P$
we have
\begin{align*}
  \chi(\P_{\ge\sigma})&=\sum_{\tau>\sigma}(-1)^{\rank\tau-1}+(-1)^{\rank\sigma-1}
  =(-1)^{\rank\sigma}\Bigl(\sum_{\tau>\sigma}(-1)^{\rank\tau-\rank\sigma-1}-1\Bigr)\\
  &=(-1)^{\rank\sigma}\Bigl(\sum_{\varnothing\ne\rho\in\lk_K\sigma}(-1)^{\dim\rho}-1\Bigr)
  =(-1)^{\rank\sigma}\bigl(\chi(\lk_K\sigma)-1\bigr)
\end{align*}
(since $K$ a simplicial complex, the poset of non-empty faces of
$\lk_K\sigma$ is isomorphic to $\P_{>\sigma}$ with shifted rank
function). Now, because $K$ is a triangulation of a closed
$(n-1)$-manifold, the link of a non-empty simplex $\sigma$ is a
homology sphere of dimension $(n-\rank\sigma-1)$. Therefore,
$\chi(\lk_K\sigma)=1+(-1)^{n-\rank\sigma-1}$ and
$\chi(\P_{\ge\sigma})=(-1)^{n-1}$ for $\sigma\ne\varnothing$. We
also have $\lk_K\varnothing=K$. It follows from
Theorem~\ref{theo:DS} that
\[
\begin{split}
\sum_{i=0}^n(h_{n-i}-h_{i})t^i
&=\big(1+(-1)^{n}(\chi(K)-1)\big)(t-1)^n\\
&=(-1)^n\big(\chi(K)-\chi(S^{n-1})\big)(t-1)^n.
\end{split}
\]
Comparing the coefficients of $t^i$ of both sides, we obtain the
corollary.
\end{proof}

For a face $F$ of a torus graph $\Gamma$, we define its Euler
number $\chi(F)$ by
\begin{equation} \label{eqn:chimF}
  \chi(F):=\sum_{H\subseteq F}(-1)^{\dim H}
\end{equation}
where $H$ is a face of $F$.

\begin{coro}
$\displaystyle{\sum_{i=0}^n(h_{n-i}-h_{i})t^i=\sum_{\mF\subseteq\mG}
\bigl(1-\chi(\mF)\bigr)(t-1)^{\dim \mF}}$.  \newline In
particular, the equations $h_i=h_{n-i}$ hold if $\chi(\mF)=1$ for
every face $\mF$ of~$\mG$.
\end{coro}
\begin{proof}
We apply Theorem~\ref{theo:DS} to the simplicial poset $\P(\mG)$
associated with the graph $\mG$. Given a face $F$, denote by
$\sigma$ the corresponding element of $\P(\mG)$. Then
$\rank\sigma=(n-\dim F)$ and
\[
  \chi(\P_{\ge\sigma})=\sum_{\tau\ge\sigma}(-1)^{\rank\tau-1}=
  \sum_{H\subseteq F}(-1)^{n-\dim H-1}=(-1)^{n-1}\chi(F).
\]
This together with Theorem~\ref{theo:DS} proves the corollary.
\end{proof}

\end{document}